\documentclass[12pt]{article}

\usepackage[english]{babel}
\usepackage{amssymb,a4wide}
\usepackage{doublespace}
\usepackage{amsmath,xypic,epsf,mathrsfs,theorem}

\usepackage[all]{xy}

\newtheorem{theorem}{Theorem}[subsection]
\newtheorem{lemma}[theorem]{Lemma}
\newtheorem{corollary}[theorem]{Corollary}

\newtheorem{proposition}[theorem]{Proposition}

{\theorembodyfont{\upshape}
\newtheorem{definition}[theorem]{Definition}
\newtheorem{remark}[theorem]{Remark}

\newtheorem{question}[theorem]{Question}}

\numberwithin{equation}{section}
\numberwithin{theorem}{section}

\newcommand{\cp}{completely positive}
\newcommand{\cpc}{completely positive contraction}

\newcommand{\sC}{{\mathscr{C}}}

\newcommand{\sK}{{\mathscr{K}}}

\newcommand{\Ker}{\mathrm{Ker}}
\newcommand{\image}{\mathrm{Im}}
\newcommand{\dist}{\mathrm{dist}}
\newcommand{\diag}{\mathrm{diag}}
\newcommand{\Id}{\mathrm{id}}

\newcommand{\spek}{\mathrm{sp}}

\newcommand{\Ad}{\mathrm{Ad}\,}

\newcommand{\cM}{{\cal M}}
\newcommand{\ep}{{\varepsilon}}

\newcommand{\calL}{{\cal L}}

\newcommand{\cO}{{\cal O}}

\newcommand{\C}{{\mathbb C}}
\newcommand{\Z}{{\mathbb Z}}
\newcommand{\N}{{\mathbb N}}

\newcommand{\cK}{{\cal K}}

\newcommand{\R}{{\mathbb R}}
\newcommand{\Cs}{{$C^*$-al\-ge\-bra}}
\newcommand{\A}{{\mathcal{A}}_{[0,1]}}
\newcommand{\sCs}{{sub-$C^*$-al\-ge\-bra}}

\newcommand{\sh}{{$^*$-ho\-mo\-mor\-phism}}

\newcommand\eqdef{{\;\overset{\mbox{\scriptsize def}}{=}\;}}

\newenvironment{proof}[1][Proof:]
{\begin{trivlist}\item[]\textbf{#1} }
{\hbox{}\nobreak\hfill\quad\hbox{$\square$}\end{trivlist}}

\begin{document}

\title{Purely infinite $C^*$-algebras: ideal-preserving zero homotopies}

\author{Eberhard Kirchberg and Mikael R\o rdam}
\date{}
\maketitle

\begin{abstract} \noindent
We show that if $A$ is a separable,
nuclear, $\cO_\infty$-absorbing (or strongly purely infinite) \Cs{}
which is \emph{homotopic to zero in an ideal-system preserving way},
then $A$ is the inductive limit of \Cs s of the form $C_0(\Gamma,v)
\otimes M_k$, where $\Gamma$ is a finite graph (and $C_0(\Gamma,v)$ is the
algebra of continuous functions on $\Gamma$ that vanish at
a distinguished point $v \in \Gamma$).

We show further that if $B$ is any
separable, nuclear \Cs{}, then $B \otimes \cO_2 \otimes \cK$ is
isomorphic to a crossed product $D \rtimes_\alpha \Z$, where $D$ is an
inductive limit of \Cs s of the form $C_0(\Gamma,v) \otimes M_k$
(and $D$ is $\cO_2$-absorbing and homotopic to zero in an ideal-system
preserving way).
\end{abstract}

\section{Introduction} \label{sec:intro}
\noindent Cuntz defined in \cite{Cuntz:KOn} a simple \Cs{} to be
purely infinite if each of its non-zero here\-ditary sub-\Cs s contain
an infinite projection. The class of simple purely infinite \Cs s has
since then received much attention; in parts because many interesting
and naturally occurring \Cs s, such as the Cuntz algebras $\cO_n$, $2
\le n \le \infty$, are purely infinite, and in parts because strong
classification theorems have been obtained for this class of \Cs s
(see \cite{Kir:fields}, \cite{Phi:class}, and \cite{Ror:encyc}).

The notion of being purely infinite was extended to non-simple \Cs s
in \cite{KirRor:pi}. A (possibly non-simple) \Cs{} $A$ is  said to be
purely infinite if $A$ has no character and if for every pair of
positive elements $a,b \in A$, such that $a$ belongs to the closed
two-sided ideal in $A$ generated by $b$, there is a sequence $\{x_n\}$
in $A$ with $x_n^*bx_n \to a$. Every purely infinite \Cs{} is
traceless (in the sense that no algebraic ideal of the \Cs{} admits a
non-zero trace or a quasitrace), and $A \otimes \cO_\infty$ is purely
infinite for every \Cs{} $A$. It was shown in \cite{KirRor:pi2} that
if $A$ is a nuclear, separable \Cs{} that is either stable or unital,
then the following three statements are equivalent: $A \cong A
\otimes \cO_\infty$, $A$ is approximately divisible and
traceless, and $A$ is \emph{strongly purely infinite} (an internal
algebraic condition, see \cite[Definition~5.1]{KirRor:pi2}). If $A$ is
of real rank 
zero, then these 
conditions are satisfied if and only if every non-zero projection in
$A$ is properly infinite.

There is a classification (in terms of an ideal related version of
Kasparov's $KK$-theory) for separable, stable, nuclear \Cs s that absorb the
Cuntz algebra $\cO_\infty$; and in the special case of
$\cO_2$-absorbing \Cs s the classifying invariant is nothing but the
primitive ideal space (or equivalently, the ideal lattice). These
results were obtained by the first named author in \cite{Kir:fields}.
(A \Cs{} $A$ is said to be $\cO_\infty$- and $\cO_2$-absorbing,
respectively, if $A \cong A \otimes \cO_\infty$ and $A \cong A \otimes
\cO_2$, respectively.) The following question is open:

\begin{question} \label{q1}
Is every point-complete\footnote{A $T_0$-space is \emph{point-complete}
   (or \emph{spectral})
   if every prime closed subset is the closure of a point. A closed
    subset $F$ of a T$_0$-space $X$ is called \emph{prime} if
   $F=G \cup H$ implies $F=G$ or $F=H$, when $G$ and $H$ are closed.},
second countable, locally quasi-compact $T_0$-space
the primitive ideal space of a separable (or separable and nuclear) \Cs{}.
\end{question}

\noindent Going a step
further, it would be desirable to have an algorithm that to each
admissible topological space
assigns a separable, nuclear, $\cO_2$-absorbing \Cs{} with that
space as its primitive ideal space.
We take a step in this direction at the end of Section 6 by showing the
following result: For every \Cs{} $B$ there is a commutative \sCs{}
$C\subseteq B\otimes \cO _2$ such that
$\Psi\colon J\mapsto C\cap (J\otimes \cO _2)$ is an
injective  lattice-morphism from the Hausdorff lattice of closed ideals
of $B$ into the lattice of closed ideals of $C$
(in particular, $\Psi (J_1+J_2)=\Psi (J_1)+\Psi (J_2)$).

In \cite{Ror:piah} a \Cs{} $\A$ was constructed that has primitive ideal space
$[0,1)$
with the non-Hausdorff $T_0$ topology
$\{ [0,\alpha) : \alpha \in
(0,1]\} \cup \{ \emptyset \}$. The \Cs{} $\A$ is
$\cO_\infty$-absorbing, and at the same time an inductive limit of \Cs
s of the form $C_0([0,1),M_{2^k})$. We observe in Section 6 that $\A$ is
zero homotopic in an ideal-system preserving way, i.e., that there is
a continuous path $\{\rho_t\}_{t \in [0,1]}$ of $^*$-endomorphisms on
$\A$ such that $\rho_0 = 0$, $\rho_1 = \Id$, and $\rho_t(I)
\subseteq I$ for every closed two-sided ideal $I$ in $\A$ and for each
$t$. By the classification of nuclear
$\cO_\infty$-absorbing \Cs s we can then conclude that $\A$ is
$\cO_2$-absorbing (see Theorem~\ref{thm:Kir}), and
hence is the unique such \Cs{}, up to isomorphism, with that
primitive ideal space.

A question we seek to answer in this paper is to what extend this
example is a special case of a more general theory. In particular we
ask the following:

\begin{question} \label{q2}
Are the following conditions equivalent for a nuclear, separable,
stable, $\cO_\infty$-absorbing \Cs{} $A$:
\begin{enumerate}
\item $A$ is zero homotopic in an ideal-system preserving way.
\item Each quotient of $A$ is homotopic to zero.
\item Each quotient of $A$ is projectionless.
\item $A$ is an AH$_0$-algebra.
\end{enumerate}
\end{question}
 An AH$_0$-algebra is a \Cs{} that is the
inductive limit of direct sums of \Cs s of the form $C_0(X,M_n)$,
where $X$ is a locally compact space. The spaces $X$, that occur in
the direct limit, can in our case be taken to be (one-dimensional)
graphs. The implications (i) $\Rightarrow$ (ii) $\Rightarrow$ (iii) and (iv)
$\Rightarrow$ (iii) are true and trivial. We show here
(Theorem~\ref{thm:pi=ah_0}) that (i)
$\Rightarrow$ (iv) if $A$ is
$\cO _\infty$-absorbing. We do not know if (iii) $\Rightarrow$ (i).

If $A$ satisfies (i) above, then so does $A \otimes B$ for every \Cs{}
$B$ (as long as one of $A$ and $B$ is nuclear). In particular, $\A
\otimes B$ satisfies (i), and hence (iv), in Question~\ref{q2} for
every separable nuclear
\Cs{} $B$. We show in Section~6 that there is an automorphism $\alpha$
on $\A$ such that $\A \rtimes_\alpha \Z$ is isomorphic to $\cO_2
\otimes \cK$. Hence, if $B$ is any separable,
nuclear \Cs{}, then $B \otimes \cO_2 \otimes \cK$ is isomorphic to $D
\rtimes_\beta \Z$, where $D = \A \otimes B$ and $\beta = \alpha \otimes
\Id_B$, and $D$ is an AH$_0$-algebra. In other words, every separable, nuclear,
stable, $\cO_2$-absorbing \Cs{} is the crossed product of an
AH$_0$-algebra by an action of the integers, thus confirming
the following fundamental question in the special case of
$\cO_2$-absorbing \Cs s.

\begin{question} \label{q3}
Is every separable, nuclear, stable \Cs{} isomorphic to a crossed
pro\-duct $D \rtimes_\alpha \Z$, where $D$ is an inductive limit of type
I \Cs s, and where $\alpha$ is an automorphism on $D$?
\end{question}
An affirmative answer to this question will imply that the Universal
Coefficient Theorem (UCT) holds for all separable nuclear \Cs s. (To
this end, it will suffice to affirm Question~\ref{q3} for \Cs s that absorb
the Cuntz
algebra $\cO_\infty$.) Unfortunately, our result for $\cO_2$-absorbing
\Cs s does not add new information to the UCT.

In Section~2 we remind the reader of some results about \cp{} maps,
and prove some sharpened versions of known results (in parts from our
earlier paper \cite{KirRor:pi2}).

In Section~3 we show that any Abelian sub-\Cs{} $B$ of an
$\cO_2$-absorbing \Cs{} $A$ is approximately contained in an Abelian
sub-\Cs{} $C$ of $A$, where the spectrum of $C$ is a
(one-dimensional) graph. We also show that if $X$ is any locally
compact metrizable space, then $C(X) \otimes 1 \subseteq D \subseteq
C(X) \otimes \cO_2$ for some Abelian \Cs{} $D$ whose spectrum is
one-dimensional.

In Section~4 we prove a Hahn-Banach type separation theorem
(implicitly contained in our earlier paper \cite{KirRor:pi2}) for
\cp{} maps from a nuclear \Cs{} $A$ into an arbitrary \Cs{} $B$: If $\sK$ is
an operator convex set of \cp{} maps from $A$ to $B$, and if $T \colon
A \to B$ is a \cp{} map, then $T$ belongs to the point-norm closure of
$\sK$ if and only if $T(a)$ belongs to the closure of $\{V(a) : V \in
\sK\}$ for all $a \in A$.  This result, and the approximation
result from Section~3, are the main technical ingredients in the proof
in Section~5 of our main result. There are also other applications of
our Hahn-Banach theorem. We give a new proof of the fundamental
uniqueness theorem for $\cO_2$-absorbing \Cs s: If $\varphi$ and
$\psi$ are \sh s from a separable, nuclear \Cs{} $A$ to an
$\cO_2$-absorbing, stable \Cs{} $B$, then $\varphi$ and $\psi$ are
approximately unitarily equivalent if and only if they induce the same
map from the ideal lattice of $B$ to the ideal lattice of $A$.

Section~5 contains the proof of our main result described in the
abstract, and Section~6 contains some applications
thereof, already mentioned.

\section{Preliminaries} \label{sec:prelim}

\noindent We remind the reader of concepts and notations
related to limit algebras.
A filter $\omega$ on $\N$ is called free if
$\bigcap_{A \in \omega} A = \emptyset$. The filter $\omega_\infty$ of
all co-finite subsets of $\N$ is free; and a filter $\omega$ is free
if and only if $\omega_\infty \subseteq \omega$.
In this paper $\omega$ always means a \emph{free} filter
on $\N$.

If $A$ is a \Cs{} and $\omega$ is a filter on $\N$, then
$A_\omega$ is defined to be the quotient-\Cs{}
$\ell_\infty(A)/I_\omega$, where $\ell_\infty(A)$ is the \Cs{} of all
bounded sequences $\{a_n\}_{n=1}^\infty$, with $a_n \in A$, and
$I_\omega$ is the closed two-sided ideal in $\ell_\infty(A)$
consisting of
those sequences $\{a_n\}_{n=1}^\infty$ for which $\lim_\omega \|a_n\|
= 0$. We let $\pi_\omega$ denote the quotient mapping $\ell_\infty(A)
\to A_\omega$.
A sequence converges along the filter $\omega_\infty$ (defined above)
if and only if it converges in the usual sense as a
sequence. Therefore $I_{\omega_\infty} = c_0(A)$, which is the set of
sequences $\{a_n\}_{n=1}^\infty$ such that $\lim_{n \to \infty}
\|a_n\| = 0$. We let $A_\infty$ denote the \Cs{}
$A_{\omega_\infty}$, i.e., $A_\infty = \ell_\infty(A)/c_0(A)$.

We shall often have occasion to consider the \Cs{} $\cM(A)_\omega$,
where $A$ is a non-unital \Cs{}. Note that there are natural inclusions
$A \subseteq A_\omega \subseteq \cM(A)_\omega\subseteq \cM(A_\omega)$.

We remind the reader of the following:

\begin{definition}[Approximately inner \cp{} maps] \label{def:app_inner}
\hfill A completely \\
positive map $T$ from a \Cs{} $A$ to another \Cs{} $B$ is called
\emph{$n$-step inner} if
  there are elements $b_1, \dots, b_n$ in $B$ such that $T(a) =
  \sum_{j=1}^n b_j^*ab_j$ for all $a \in A$. $T$ is called
  \emph{inner} if it is $n$-step
inner for some $n$. If $T$ is the point-norm limit of a sequence of
$n$-step inner or inner \cp{} maps, then $T$ is called
\emph{approximately $n$-step inner} and \emph{approximately
  inner}, respectively.
\end{definition}

\begin{lemma}[cf.\ Lemma 7.4 of \cite{KirRor:pi2}] \label{lm:1-step-inner}
Let $A$ be a stable \Cs{}, let $\omega$ be a filter on $\N$, and
let $B$ be a separable sub-\Cs{} of $A_\omega$. Then each
approximately 1-step inner \cpc{} $T \colon B \to A_\omega$ is of the
form $T(b) = t^*bt$, $b \in B$, for some isometry $t \in \cM(A)_\omega$.
\end{lemma}

\begin{proof} By assumption and \cite[Lemma~7.2]{KirRor:pi2} there is
  a sequence of \emph{contractions} $\{f_n\}$ in $A_\omega$ such that
  $f_n^*bf_n \to T(b)$ for all $b$ in $B$.
  By \cite[Lemma~2.5]{KirRor:pi2} there is a single contraction $d$ in
  $A_\omega$ for which we have the exact relation $T(b) = d^*bd$ for
  all $b$ in $B$.
  Write $d = \pi_\omega(d_1,d_2, \dots)$, where
  each $d_n$ is a contraction in $A$. Let $\{b_1,b_2, \dots\}$ be a
  countable dense subset of $B$ and write $b_j =
  \pi_\omega(b_{j,1},b_{j,2}, \dots)$, where $b_{j,n} \in A$. As in
  the proof of \cite[Lemma~7.4]{KirRor:pi2}, using stability of $A$,
  we can find isometries $s_{n,1}$ and $s_{n,2}$ in $\cM(A)$ for each
  $n \in \N$  satisfying $s_{n,1}s_{n,1}^* + s_{n,2}s_{n,2}^*=1$ and
  $\|b_{j,n}s_{n,2}s_{n,2}^*\|  \le 1/n$ for
    $j = 1,2, \dots, n$. As in the proof of \cite[Lemma~7.4]{KirRor:pi2}, set
$$
t_n = s_{n,1}s_{n,1}^*d_n + s_{n,2} \big(1 -
d_n^*s_{n,1}s_{n,1}^*d_n\big)^{1/2} \in \cM(A),  \qquad
t = \pi_\omega(t_1,t_2, \dots) \in \cM(A)_\omega.
$$
Each $t_n$ is an isometry and
$\|t_n^*b_{j,n}t_n - d_n^*b_{j,n}d_n\| \to 0$ as $n \to \infty$ for
all $j$. It follows that $t$ is an
isometry and that
$$t^*b_jt = \pi_\omega(t_1^*b_{j,1}t_1,t_2^*b_{j,2}t_2, \dots) =
\pi_\omega(d_1^*b_{j,1}d_1,d_2^*b_{j,2}d_2, \dots) = d^*b_jd =
T(b_j)$$
for all $j$, and so $T(b)=t^*bt$ for all $b$ in $B$.
\end{proof}

\begin{lemma} \label{lm:O_2-in-commutant}
Let $A$ be a \Cs{}, let $\omega$ be a free filter on $\N$, and let $B$
be a separable sub-\Cs{} of $A_\omega$.
\begin{enumerate}
\item If $A \cong A \otimes \cO_2$, then there is a unital embedding
of $\cO_2$ into $B' \cap \cM (A)_\omega$.
\item If $A \cong A \otimes \cO_\infty$, then there is a unital embedding
of $\cO_\infty$ into $B' \cap \cM (A)_\omega$.
\end{enumerate}
\end{lemma}

\begin{proof}
We prove (i) and (ii) simultaneously. Let $m$ represent either $2$ or
$\infty$ and assume that $A \cong A \otimes \cO_m$. Then
there is a sequence $\varphi_n \colon \cO_m \to \cM(A)$ of unital \sh
s that asymptotically commutes with $A$, ie.,
$\varphi_n(x)a-a\varphi_n(x) \to 0$ for all $x \in \cO_m$ and all $a
\in A$ (see eg.\ \cite[Theorem 7.2.2 and Remark 7.2.3]{Ror:encyc}).
Let $\pi_\omega \colon \ell_\infty(A) \to A_\omega$ be the quotient
mapping, let $\{b_1,b_2,b_3, \dots \}$ be a dense subset of $B$, let
$\{x_1,x_2,x_3, \dots\}$ be a dense subset of $\cO_m$, and write
$b_j = \pi_\omega(b_{j1},b_{j2},b_{j3}, \dots)$, where $b_{ji} \in
A$.
For each natural number $k$ find a natural number $n_k$ such that
$\|\varphi_{n_k}(x_i)b_{jk}-b_{jk}\varphi_{n_k}(x_i)\| < 1/k$ for
$i, j=1,2, \dots, k$. Define $\varphi \colon \cO_m \to
\cM(A)_\omega$ by
$$\varphi(x) =
\pi_\omega(\varphi_{n_1}(x),\varphi_{n_2}(x),\varphi_{n_3}(x), \dots),
\qquad x \in \cO_m.$$
Because $\omega$ is free and
$$\lim_{k \to \infty} \|\varphi_{n_k}(x_i)b_{jk} -
b_{jk}\varphi_{n_k}(x_i)\| = 0, \qquad i,j \in \N,$$
it follows that $\varphi(x_i)b_j=b_j\varphi(x_i)$ for all $i$ and $j$;
and this shows that the image of $\varphi$ commutes with $B$.
\end{proof}

\noindent A \Cs{} $A$ is called \emph{$\cO_2$-absorbing} or
\emph{$\cO_\infty$-absorbing}, if $A \cong A \otimes
\cO_2$ or $A \cong A \otimes \cO_\infty$, respectively. We have $\cO_2
\cong \cO_2 \otimes \cO_2 \otimes \cO_\infty$ and $\cO_\infty \cong
\cO_\infty \otimes \cO_\infty$ (see \cite[Theorem~3.8 and
Theorem~3.15]{KirPhi:classI}), so $A \otimes \cO_2$ is both $\cO_2$-
and $\cO_\infty$-absorbing, and $A \otimes \cO_\infty$ is
$\cO_\infty$-absorbing for every \Cs{} $A$

Let $D$ be a unital \Cs{} and suppose that
$v_1, v_2$ are isometries in $D$ satisfying the $\cO_2$-relation:
$v_1v_1^*+v_2v_2^*=1$. We shall then consider the \emph{Cuntz sum} of elements
$a,b \in D$ defined by $a \oplus_{v_1,v_2} b = v_1av_1^* + v_2bv_2^*$.

\begin{lemma} \label{lm:Cuntz-sum}
Let $D$ be a unital \Cs{}, let $s$ be an isometry in $D$, and let
$V(a)=s^*as$ be the corresponding unital completely positive map on $D$.
Suppose that
$v_1, v_2$ are isometries in $D$ satisfying the $\cO_2$-relation.
Put $w_1 = (1-ss^*) + sv_1s^*$ and $w_2 = sv_2$. 

Then $w_1,w_2$ are isometries satisfying
the $\cO_2$-relation, and 
\begin{eqnarray*}
& \|a \oplus_{w_1,w_2} V(a) - a\|  \le  \|[v_1,V(a)]\| + \|[v_2,V(a)]\|
+ 2\|[a,ss^*]\|, & \\
& \|[a,ss^*]\|  =  \max\{\|V(a^*a)-V(a)^*V(a)\|^{1/2},
\|V(aa^*)-V(a)V(a)^*\|^{1/2}\},&
\end{eqnarray*}
for all $a \in D$.
\end{lemma}

\begin{proof}
It is straightforward to check that $w_1,w_2$ satisfy the
$\cO_2$-relation. Next,
\begin{eqnarray*}
 \|a \oplus_{w_1,w_2} V(a) -a\| 
& \le & \|(1-ss^*)a(1-ss^*) +
s\big(s^*as \oplus_{v_1,v_2} s^*as\big)s^*-a\| \\ &&
 + \; \|(1-ss^*)asv_1^*s^*+sv_1s^*a(1-ss^*)\| \\ 
& \le & \|(1-ss^*)a(1-ss^*) + ss^*ass^*-a\| \\
&&+ \; \|V(a) \oplus_{v_1,v_2}
V(a) - V(a)\| +  \|[a,ss^*]\| \\
& \le & \|[v_1,V(a)]\|
+\|[v_2,V(a)]\| +  2\|[a,ss^*]\|.
\end{eqnarray*}
The last statement follows from the identities
\begin{eqnarray*}
\|[a,ss^*]\| & = & \max\{ \|(1-ss^*)ass^*\|, \|ss^*a(1-ss^*)\|\},\\
\|V(a^*a)-V(a)^*V(a)\|  & = & \|s^*a^*(1-ss^*)as\|  \; = \;
\|(1-ss^*)as\|^2  \; = \; \|(1-ss^*)ass^*\|^2,\\
\|V(aa^*)-V(a)V(a)^*\| & = &  \|(1-ss^*)a^*ss^*\|^2 \;= \;
\|ss^*a(1-ss^*)\|^2.
\end{eqnarray*}
\end{proof}

\begin{lemma}[cf.\ Lemma~1.12 of \cite{KirPhi:classI}]
\label{lm:cantor-bernstein}
Suppose that $A$ is a separable, stable,
$\cO_2$-absorbing \Cs{} and that
$s$ and $t$ are isometries in $\cM(A)$.
Let $V(a)=s^*as$ and $W(a)=t^*at$ be
the corresponding unital completely positive
maps on $\cM (A)$.
Then for every 
finite subset $F$ of $A$ and
for every $\varepsilon >0$
there is a unitary element $u$ in
$\cM(A)$ with $$\|u^*V (a)u - a\| \le 5\kappa +\ep$$
for all $a \in F$, where 
$$\kappa = \max_{a \in F \cup F^*} \{ \|V(a^*a)-V(a^*)V(a) \|^{1/2},
\|WV(a^*a) - WV(a)^*WV(a)\|^{1/2},\| WV(a)-a\|\}.$$
\end{lemma}

\begin{proof} 
Find  isometries $r_1,r_2$ in $\cM(A)$ satisfying the $\cO_2$-relation
such that 
\begin{equation} \label{eq:O_2}
\| [r_1,V(a)] \| + \|[r_2,V(a)]\| \le \ep/2, \qquad  \|[r_1,WV(a)]\| +
\|[r_2,WV(a)]\| \le \ep/2,
\end{equation}
for all $a \in F$ (see eg.\ \cite[Theorem 7.2.2 and Remark
7.2.3]{Ror:encyc}).
Define new sets of isometries $(t_1,t_2)$ and
$(s_1,s_2)$ satisfying the $\cO_2$-relation by 
$$t_1=(1-tt^*)+tr_1t^*, \quad t_2=tr_2, \qquad s_1=(1-ss^*)+sr_1s^*,
\quad s_2=sr_2,$$
(cf.\ Lemma~\ref{lm:Cuntz-sum}),
and set $u= t_1s_2^*+t_2s_1^*$. Then $u$ is unitary element in
$\cM(A)$, and $u^*(a\oplus _{t_1,t_2}b)u=b\oplus_{s_1,s_2} a$ for all $a,b
\in \cM(A)$. Thus 
\begin{eqnarray*}
\| u^*V(a)u -a\|  & = & \| u^*V(a)u -u^* \big(V(a)\oplus
_{t_1,t_2}WV(a)\big)u + WV(a) \oplus_{s_1,s_2} V(a)  -a \| \\
& \le & \|V(a) - V(a)\oplus_{t_1,t_2}WV(a)\| + 
\|a \oplus_{s_1,s_2} V(a) -a\| + \|WV(a)-a\|
\end{eqnarray*}
for $a\in \cM(A)$. By the inequality
$$0 \le W(V(a)^*V(a)) - WV(a)^*WV(a) \le WV(a^*a)-WV(a)^*WV(a),$$ 
and by Lemma~\ref{lm:Cuntz-sum} and \eqref{eq:O_2}
we see that the two first terms above are each bounded by $2\kappa + \ep/2$,
when $a \in F$, and the last term is at most $\kappa$. This proves the
lemma.
\end{proof}

\begin{remark}\label{rem:cantor-bernstein}
  Lemma~\ref{lm:cantor-bernstein} holds with $\cM(A)$ replaced with
  $\cM(A)_\omega$ for any free filter $\omega$, even with
  $\ep=0$. Indeed, \eqref{eq:O_2} is satisfied with $\ep=0$ by
  Lemma~\ref{lm:O_2-in-commutant}~(i) (when $B$ is the separable
  sub-\Cs{} of $A_\omega$ generated by $V(A) \cup WV(A)$). The rest of
  the proof of Lemma~\ref{lm:cantor-bernstein} works with $\cM(A)$
  replaced with $\cM(A)_\omega$. 

It follows in particular that if $\{s_n\}$ and $\{t_n\}$ are
isometries in $\cM(A)_\omega$ such that
$$\lim_{n \to \infty} \|s_n^*abs_n-s_n^*as_ns_n^*bs_n\| = 0 , \qquad \lim_{n
  \to \infty} \|t_n^*s_n^*as_nt_n-a\|=0 \quad \text{for all} \; a,b
\in A,$$
then there is a sequence $\{u_n\}$ of unitaries
in $\cM(A)_\omega$ 
such that $\|u_n^*s_n^*as_nu_n -a\| \to 0$ for all $a \in A$.
\end{remark}

\section{Local and global AH$_0$-algebras} \label{sec:ah_0}

\noindent The main result of this section is
Theorem~\ref{thm:ah_0} which states that any $\cO_2$-absorbing \Cs{},
that locally is an AH$_0$, is actually an inductive limit of \Cs s of
the form $C_0(\Gamma,v) \otimes M_k$, where $(\Gamma,v)$ is a finite
pointed graph. Along the way we prove a perturbation result, which
states that a
subalgebra $B$ of an $\cO_2$-absorbing \Cs{}, with $B \cong C_0(X,M_k)$,
is approximately contained in a subalgebra that is isomorphic to
$C_0(Y,M_k)$ for some one-dimensional space $Y$. We end the section by
showing that for every compact metrizable space $X$ there is a
one-dimensional compact metrizable space $Y$ such that $C(X) \otimes 1
\subseteq C(Y) \subseteq C(X) \otimes \cO_2$.

\begin{lemma} \label{lm:normal}
Let $X$ be a compact subset of the complex plane and let $f \in
C(X,\cO_2)$ be the function $f(z) = z1$. Then for each $\lambda \in \C$
and for each $\ep >0$ there is a normal element $g \in C(X,\cO_2)$
that satisfies $\lambda \notin \spek(g)$, $\|f-g\| \le 2\ep$, and $g(z)
= f(z)$ when $|z-\lambda| \ge \ep$.
\end{lemma}

\begin{proof} Let $D=\{ z\in \C: |z-\lambda| \leq \ep \}$ be the
closed disk with radius $\ep$ and center $\lambda$, and let $\partial
D$ denote its boundary. The map $u_0 \in C(\partial D, \cO _2)$ given
by $u_0(z) = (z-\lambda)/\ep$ is unitary. The restriction
map $C(D,\cO _2)\to C(\partial D,\cO _2)$ is a $^*$-epimorphism,
and the unitary group of $C(\partial D, \cO _2)$ is connected
(cf.~\cite{Cuntz:KOn}), so there is an element
$u \in C(D,\cO _2)$ with values in the unitary group
of $\cO _2$ such that $u(z)=u_0(z)$ for $z\in \partial D$.
The map $g\colon X\to \cO _2$ given by
$$g(z) = \begin{cases}  \ep u(z) +\lambda 1, & z\in X\cap D, \\
z1, & z\in X\setminus D, \end{cases}$$
defines a normal element in $C(X,\cO _2)$
(because $g(z)$ is normal for every $z\in X$). Use that
$g(z) - \lambda 1$ is invertible and $\|g(z)-f(z)\| = \|g(z)-z1\|
\le 2\ep$ for every $z \in D$ to see that $g$ has the desired properties.
\end{proof}

\noindent Let $a_1,a_2, \dots, a_n$ be commuting normal elements in a
\Cs{} $A$. Their joint spectrum, $\spek(a_1,a_2,
\dots, a_n)$, is the compact subset of $\C^n$ of all $n$-tuples
$(\rho(a_1), \rho(a_2), \dots, \rho(a_n))$, where $\rho$ is a character
on $C^*(a_1,a_2, \dots, a_n)$. Each continuous function $f \colon
\spek(a_1,a_2, \dots, a_n) \to \C$ defines an element
$f(a_1,a_2, \dots, a_n) \in C^*(a_1,a_2,\dots,a_n)$ that satisfies
$\rho(f(a_1,a_2, \dots,a_n)) = f(\rho(a_1), \rho(a_2), \dots,
\rho(a_n))$ for all characters $\rho$ on $C^*(a_1,a_2,\dots,a_n)$.

\begin{lemma} \label{lm:normal2}
Let $A$ be a \Cs{}, let $a_1,a_2$ be two commuting self-adjoint
elements in $A$, and suppose
there is a unital \sh{} $\varphi \colon \cO_2 \to \cM(A) \cap
\{a_1,a_2\}'$. Then for each point $(s,t) \in \R^2 \setminus
\{(0,0)\}$ and for each
$\ep >0$ there are self-adjoint elements $b_1,b_2$ in $A$ such that
$(s,t) \notin \spek(b_1,b_2)$ and $\|a_j-b_j\| < \ep$ for $j=1,2$.
\end{lemma}

\begin{proof} Put $a = a_1+ia_2$ and put $\lambda = s+it \ne 0$. The \Cs{}
  $C^*(a, \varphi(\cO_2))
  \subseteq \cM(A)$ is isomorphic to $C(X,\cO_2)$,
  where $X = \spek(a)$, via an isomorphism that maps $a$ to the
  function $f$ given by $f(z)=z1$. Take $\eta >0$ such that $\eta <
  \ep/2$ and $\eta \le |\lambda|$. Use Lemma~\ref{lm:normal} to find a
normal element $g$ in $C(X,\cO_2)$ with $\lambda \notin
\spek(g)$, $\|f-g\| \le 2\eta$, and $g(z)
= f(z)$ when $|z-\lambda| \ge \eta$. Let $b$ in $C^*(a,
\varphi(\cO_2))$ correspond to $g$, so that $\|a-b\|
  \le 2\eta < \ep$ and $\lambda \notin \spek(b)$. Since $g(0)=f(0)=0$,
  because $|0-\lambda| \ge
  \eta$, we conclude that $g$ belongs to the ideal in
  $C(X,\cO_2)$ generated by $f$. Thus $b$ belongs to the
  ideal in $C^*(a, \varphi(\cO_2))$ generated by $a$, and this ideal
  is contained in $A$, whence $b$ belongs to $A$.

Write $b = b_1 + ib_2$ where $b_1$ and $b_2$ are self-adjoint. Then $b_1$
commutes with $b_2$ (because $b$ is normal), $\|a_j-b_j\| \le
\|a-b\|<\ep$, and $(s,t) \notin \spek(b_1,b_2)$---the latter because
$(s,t) \in \spek(b_1,b_2)$ if and only if $s+it \in \spek(b_1+ib_2)$.
\end{proof}

\noindent
For each $\ep >0$ and for each natural number $n$, let $\Gamma_{\ep,n}$
denote the 1-dimensional grid in $\R^n$ consisting of those $n$-tuples
$(t_1,t_2,
\dots, t_n)$ for which $t_j$ belongs to $\ep \Z$ for all but at most one
$j$. Every connected compact subset of $\Gamma_{\ep,n}$ is
(homeomorphic to) a finite graph. 

A finite graph is (here) a topological space
that consists of finitely many vertices, each homeomorphic to a point,
and finitely many edges, each homeomorphic to an interval, so that each
edge connects two vertices. There are no crossings of edges (except at
vertices). Any finite graph $\Gamma$ is a compact Hausdorff space. In
this paper we shall exclusively be concerned with finite graphs, and
any graph will tacitly be understood to be finite.
 
If $v$ is a point in $\Gamma$ (a vertex or
not), then $C_0(\Gamma, v)$ denotes the set of continuous functions
on $\Gamma$ that vanish at $v$. Equivalently, $C_0(\Gamma, v) =
C_0(\Gamma \setminus \{v\})$. The pair $(\Gamma, v)$ is called a
\emph{pointed graph}.

We state a topological lemma which is needed for
Lemma~\ref{lm:holes} below. For each $t=(t_1,t_2, \dots, t_n) \in \R^n$
let $\nu_n(t)$ be the number 
of $j$'s for which $t_j \in (\Z + 1/2)\ep$. Let $X_{\ep,n}$ be the set
of all $t \in \R^n$ for which $\nu_n(t) \le 1$. Observe that
$\Gamma_{\ep,n} \subseteq X_{\ep,n} \subseteq \R^n$.

\begin{lemma} \label{lm:top} For each natural number $n \ge 2$ and for each
  $\ep >0$ there is a continuous retract $F \colon X_{\ep,n} \to
  \Gamma_{\ep,n}$ that satisfies $\|F(t)-t\|_\infty \le \ep$ for all
  $t \in X_{\ep,n}$.
\end{lemma}

\begin{proof}
For $1 \le k \le n$, let $Z_{\ep,n,k}$ be the set of all $t =(t_1,
\dots, t_n) \in \R^n$ such that the number of $j$'s for which $t_j$
belongs to $\ep \Z$ is at least $n-k$. Put
$Y_{\ep,n,k} = Z_{\ep,n,k} \cap X_{\ep,n}$. Then 
$$X_{\ep,n} = Y_{\ep,n,n} \supseteq Y_{\ep,n,n-1} \supseteq \cdots
\supseteq Y_{\ep,n,2} \supseteq Y_{\ep,n,1} = Z_{\ep,n,1} =
\Gamma_{\ep,n}.$$
We define the function $F$ to be the composition of
continuous retracts
$$\xymatrix@C+.7pc{Y_{\ep,n,n} \ar[r]^-{F_n} & Y_{\ep,n,n-1}
  \ar[r]^-{F_{n-1}} & Y_{\ep,n,n-2} \ar[r]^-{F_{n-2}} &
  \cdots \ar[r] & Y_{\ep,n,2} \ar[r]^-{F_{2}} & Y_{\ep,n,1},}$$
that are to be constructed. It will map each cube 
$\prod_{j=1}^n [m_j\ep,(m_j+1)\ep]
\cap X_{\ep,n}$ (where $m_j$ are integers) into
itself, and any function $F$ 
with this property will satisfy
$\|F(t)-t\|_\infty \le \ep$.

We construct first a continuous retract $G_k \colon X_{\ep,k} \to
Y_{\ep,k,k-1}$ for each $k=2,3, \dots,n$. Take $t \in X_{\ep,k}$
and find integers $m_j$ such that $t$ belongs to the cube
$C=\prod_{j=1}^{k} [m_j\ep,(m_j+1)\ep]$. Observe that $X_{\ep,k}
\cap \partial C = Y_{\ep,k,k-1} \cap C$, where $\partial C$ is the
boundary of $C$. Put
$$q_C = \big((m_1+1/2)\ep, (m_2 + 1/2)\ep, \dots, (m_{k}+1/2)\ep\big)
\in C,$$
and let $\ell_{t,q_C}$ be the half-line in $\R^{k}$ that starts at
$q_C$ and runs through $t$ (note that $t \ne q_C$ because $\nu_{k}(t) \le
1$ and $\nu_{k}(q_C) = k >1$). The half-line $\ell_{t,q_C}$ intersects
$\partial C$ in exactly one point, and this
point we define to be $G_k(t)$. The point $G_k(t)$ belongs to
$X_{\ep,n}$ because $\nu_{k}(G_k(t)) =
\nu_{k}(t) \le 1$.
Note that $G_k(t)=t$ if and only if $t
\in Y_{\ep,k,k-1}$. (For $t \in C \cap X_{\ep,k}$ this happens
if and only if $t \in \partial C$.) 
This shows that $G_k$ is a continuous retract from
$X_{\ep,k}$ onto $Y_{\ep,k,k-1}$.

Put $F_n = G_n$. Let now $2 \le k \le n-1$ be given. 
The space $Z_{\ep,n,k} \supset Y_{\ep,n,k}$ is the (non-disjoint)
union of infinitely many subsets each of which is homeomorphic to
$\R^{k}$. We define below $F_k$ so that it corresponds
to the retract $G_k$ on each of
these subsets. In more detail, for each subset $\alpha$ of $\{1,2,
\dots, n\}$ with $|\alpha| = k$ let $\pi_\alpha \colon \R^n \to \R^k$
and $\pi'_\alpha \colon \R^n \to \R^{n-k}$ be the projection maps onto
the coordinates corresponding to $\alpha$ and to the complement of
$\alpha$, respectively. Then $Z_{\ep,n,k}$ is the union of the
collection of sets $D_{\alpha,m} = 
(\pi'_\alpha)^{-1}(m)$, with $\alpha$ as above and $m \in (\ep
\Z)^{n-k}$; and $\pi_\alpha|_{D_{\alpha,m}}
\colon D_{\alpha,m} \to \R^k$ is a
homeomorphism. As $\pi_\alpha(D_{\alpha,m} \cap
X_{\ep,n}) = X_{\ep,k}$ and $\pi_\alpha(D_{\alpha,m} \cap Y_{\ep,n,k-1}) =
Y_{\ep,k,k-1}$, there is a unique continuous retract $H_{\alpha,m}$
making the diagram
$$\xymatrix@C+1pc{D_{\alpha,m} \cap X_{\ep,n} \ar[d]_-{\pi_\alpha}
  \ar[r]^-{H_{\alpha,m}} & D_{\alpha,m} \cap Y_{\ep,n,k-1} 
\ar[d]^-{\pi_\alpha} \\  X_{\ep,k} \ar[r]_-{G_k} & Y_{\ep,k,k-1}}$$
commutative. Note that $H_{\alpha,m}(D_{\alpha,m} \cap C) \subseteq C$
whenever $C$ is a cube $\prod_{j=1}^n [k_j\ep, (k_j+1)\ep]$ with $k_j
\in \Z$. If $(\alpha,m) \ne (\alpha',m')$, then $D_{\alpha,m} \cap
D_{\alpha',m'}$ is (empty or) contained in 
$Y_{\ep,n,k-1}$, whence $H_{\alpha,m}(t) = H_{\alpha',m'}(t) = t$ for
$t$ in $D_{\alpha,m} \cap D_{\alpha',m'}$. 
The maps $H_{\alpha,m}$ therefore extend to a
continuous retract $F_k$ from $Y_{\ep,n,k}= 
\bigcup_{\alpha,m} D_{\alpha,m} \cap X_{\ep,n}$ onto $Y_{\ep,n,k-1} =
\bigcup_{\alpha,m} D_{\alpha,m} \cap Y_{\ep,n,k-1}$. 
\end{proof}

\begin{lemma} \label{lm:holes}
Let $a_1, a_2, \dots, a_n$ be commuting self-adjoint elements in a
\Cs{} $A$, and let $\ep >0$. Suppose that
$$\spek(a_i,a_j) \cap (\ep/2 + \ep \Z) \times (\ep/2 + \ep \Z) =
\emptyset$$
for all pairs $(i,j)$ with $i \ne j$. Then there are commuting
self-adjoint elements $b_1,b_2, \dots, b_n$ in $C^*(a_1,a_2, \dots,
a_n) \subseteq A$ such that $\|a_j-b_j\| \le \ep$ for all $j$ and
$\spek(b_1,b_2, \dots, b_n) \subseteq \Gamma_{\ep,n}$, where
$\Gamma_{\ep,n}$ is the one-dimensional grid in $\R^n$ defined above
Lemma~\ref{lm:top}. 
\end{lemma}

\begin{proof} If $n=1$, then we can take $b_1=a_1$ (because
  $\Gamma_{\ep,1} = \R$). We assume below
  that $n \ge 2$. The condition on the pairwise joint spectra of the elements
  $a_1,a_2, \dots, a_n$ ensures that the joint spectrum $\spek(a_1,a_2,
  \dots, a_n)$ is contained in the set $X_{\ep,n}$ defined above
  Lemma~\ref{lm:top}. Let $F \colon X_{\ep,n} \to \Gamma_{\ep,n}$, with
  $\|F(t)-t\|_\infty \le \ep$ for all $t \in X_{\ep,n}$, be
  the continuous retract found in Lemma~\ref{lm:top}. 
   Write $F(t) = (f_1(t),f_2(t),
  \dots, f_n(t))$, where each $f_j \colon X_{\ep,n} \to \R$ is a continuous
  function. Let $p_j \colon \R^n \to \R$ be the $j$th
  coordinate function so that
  $a_j = p_j(a_1,a_2, \dots,a_n)$ for $j=1,2,\dots,n$. Put $b_j =
  f_j(a_1,a_2, \dots, a_n)$. Then $b_1,b_2, \dots, b_n$ are
  self-adjoint (and necessarily commuting) elements in $C^*(a_1,a_2,
  \dots, a_n)$,
$$\|a_j-b_j\| = \sup\{|p_j(t)-f_j(t)| : t \in \spek(a_1,a_2, \dots,
a_n) \} \le \sup\{\|F(t)-t\|_\infty : t \in X_{\ep,n}\} \le \ep,$$
and
$$\spek(b_1,b_2, \dots, b_n) = F(\spek(a_1,a_2, \dots,a_n)) \subseteq
F(X_{\ep,n}) \subseteq \Gamma_{\ep,n}$$
as desired.
\end{proof}

\begin{proposition} \label{prop:perturbation} Let $A$ be a separable
  \Cs{} that absorbs $\cO_2$, and let
  $\omega$ be a free filter on $\N$. Suppose that $B$ is a sub-\Cs{} of
  $A_\omega$ and that $B$ is isomorphic to $C_0(X,M_k)$ for some
  locally compact 
  Hausdorff space $X$ and for some natural number $k$. Let $b_1,b_2,
  \dots, b_n$ be elements in $B$ and let $\ep>0$.

Then there exist a
  sub-\Cs{} $B_1$ of $A_\omega$ and elements $c_1,c_2, \dots, c_n \in
  B_1$ such that $\|b_j-c_j\| < \ep$ for all $j$, and such that $B_1$ is 
  isomorphic to $C_0((\Gamma,v), M_k)$, or to $C(\Gamma,M_k)$, for some
  compact subset
  $\Gamma$ of $\Gamma_{\ep,m}$, for some $v \in \Gamma$, and for some
  natural numbers $m, k$.
\end{proposition}

\begin{proof} 
Use Lemma~\ref{lm:O_2-in-commutant} to find a
  unital \sh{} $\varphi \colon \cO_2 \to \cM(A)_\omega \cap B'$.
Let $I$ be the closed
  two-sided ideal in
  $C^*(B,\varphi(\cO_2))
  \subseteq \cM(A)_\omega$
  generated by $B$. Then
$$I \cong M_k \otimes C_0(X) \otimes
  \cO_2 \quad \text{and} \quad I \subseteq A_\omega.$$
The elements $b_1,b_2, \dots, b_n$ in $B$ correspond to
  elements $\tilde{b}_1,\tilde{b}_2, \dots, \tilde{b}_n$ in $M_k
  \otimes C_0(X) \otimes 1$. Find self-adjoint elements $f_1,f_2,
  \dots, f_m$ in $C_0(X)$ such that
$M_k \otimes C^*(f_1,f_2, \dots,f_m) \otimes 1$ contains the elements
$\tilde{b}_1, \tilde{b}_2, \dots, \tilde{b}_n$. We
find pairwise commuting self-adjoint elements $g_1,g_2, \dots, g_m$ in
  $C_0(X) \otimes
  \cO_2$ with $\|f_j \otimes 1 -g_j\| < \delta$, for $j=1,2, \dots, m$
and for some small enough $\delta >0$, such that
$(\Gamma =) \; \spek(g_1,g_2, \dots,g_m) \subseteq
  \Gamma_{\ep,m}$. The \Cs{} $C^*(g_1, \dots, g_m)$ is
isomorphic to
  $C_0(\Gamma,v)$ if the origin $v = (0,\dots,0)$ belongs to
  $\Gamma$, and it is isomorphic to $C(\Gamma)$ if $v$ does not
  belong to $\Gamma$. In both cases, the \Cs{} $B_1 \subseteq 
A_\omega$ that corresponds to $M_k \otimes C^*(g_1,g_2, \dots, g_m)
  \subseteq M_k \otimes C_0(X) \otimes \cO_2$ will be as desired.

To find  $g_1, g_2,\dots, g_m$ it suffices by Lemma~\ref{lm:holes} 
to find self-adjoint elements $h_1,h_2, \dots, h_m$ in 
$C_0(X) \otimes \cO_2$ such that $\|f_j
\otimes 1 -h_j\| < \ep/2$ for all $j$ and such that $\spek(h_i,h_j)
\cap \Lambda = \emptyset$ for $i
\ne j$, where $\Lambda =(1/2 + \Z)\frac{\ep}{2} \times 
(1/2 + \Z)\frac{\ep}{2}$.

Let $R$ be the set of all commuting
$m$-tuples of self-adjoint elements in $C_0(X) \otimes \cO_2$, and 
for $i,j = 1, \dots, m$, with $i \ne j$, and for $\lambda \in \Lambda$ let
$R_{i,j,\lambda}$ be the set of all $(h_1,
\dots, h_m)$ in $R$ such that
$\lambda \notin \spek(h_i,h_j)$. Then
$R$ is a closed subset of $\big(C_0(X) \otimes \cO_2\big)^m$, whence
$R$ is a complete metric space. 
Each $R_{i,j,\lambda}$ is clearly open in $R$. We show below that
$R_{i,j,\lambda}$ is also dense in $R$. 

There is a sequence $\{\psi_\ell\}_{\ell =1}^\infty$ of
$^*$-endomorphisms on $\cO_2$ such that $\cO_2
\cap \image(\psi_\ell)'$ contains a unital copy of $\cO_2$ for each $\ell$,
and such that 
$\psi_\ell(x) \to x$ for all $x \in \cO_2$. (Indeed, since  $\cO_2 \cong
\bigotimes_{r=1}^\infty \cO_2$ (cf.\
\cite[Corollary~5.2.4]{Ror:encyc}) we can take $\psi_\ell$ to be the
endomorphism 
that corresponds to the endomorphism on $\bigotimes_{r=1}^\infty \cO_2$
that fixes the first $\ell$ copies of $\cO_2$ and shifts the remaining
copies one place to the right.) 
Put $\varphi_\ell = \Id_{C_0(X)} \otimes \psi_\ell$. Then there is a unital
copy of $\cO_2$ in $\cM(C_0(X) \otimes \cO_2) \cap \image(\varphi_\ell)'$
for every $\ell$. 

Let $(d_1, \dots, d_m) \in R$ and $\ep>0$ be given. Choose $\ell$ such that
$(\varphi_\ell(d_1), \dots, \varphi_\ell(d_m)) \in R$ is within
distance $\ep/2$ from $(d_1, \dots, d_m)$. By
Lemma~\ref{lm:normal2} there is a  commuting pair $(e_i, e_j)$ 
of self-adjoint elements in $C_0(X) \otimes \cO_2 \cap \image(\varphi_\ell)'$
within distance $\ep/2$ from $(\varphi_\ell(d_i),\varphi_\ell(d_j))$ 
and with $\lambda \notin
\spek(e_i,e_j)$. Put $e_k = \varphi_\ell(d_k)$ if $k \ne i,j$. Then $(e_1,
\dots, e_m)$ belongs to $R_{i,j,\lambda}$ and has
distance at most $\ep$ to $(d_1, \dots, d_m)$. 

Now, $R$ has the Baire property, being a complete metric space, so 
the intersection $R_0$ of all $R_{i,j,\lambda}$, where $i,j = 1, \dots, m$,
$i \ne j$, and $\lambda \in \Lambda$, is dense in $R$. This ensures
the existence of $(h_1, \dots, h_m) \in R_0$ such that $\|f_j
\otimes 1 -h_j\| < \ep/2$ for all $j$, thus completing the proof.
\end{proof}

\begin{theorem} \label{thm:ah_0}
Let $A$ be a separable \Cs{} and let $A_\infty$ be the limit algebra
$\ell_\infty(A)/c_0(A)$. The following conditions are equivalent:
\begin{enumerate}
\item $A$ is the inductive limit of a sequence
$$\xymatrix{C_0((\Gamma_1,v_1),M_{k_1}) \ar[r]^{\varphi_1} &
  C_0((\Gamma_2,v_2),M_{k_2}) \ar[r]^{\varphi_2}  &
  C_0((\Gamma_3,v_3),M_{k_3})
      \ar[r]^-{\varphi_3} &  \cdots \ar[r] & A,}$$
where each $(\Gamma_j,v_j)$ is a pointed graph, $k_j$ is a
natural number, and each connecting map $\varphi_j$ is a \sh.
\item For every finite set $a_1, a_2, \dots, a_n$ in $A$ and for every
  $\ep >0$ there exist a \Cs{} $B \cong C_0((\Gamma,v),M_k)$, for some
  pointed graph $(\Gamma,v)$ and for some $k \in \N$, a \sh{} $\varphi
  \colon B \to A$, and elements $b_1,b_2,
  \dots, b_n$ in $B$, such that $\|a_j-\varphi(b_j)\| < \ep$ for all $j$.
\item For every finite set $a_1, a_2, \dots, a_n$ in $A$ and for every
  $\ep >0$ there exist a \Cs{} $B \cong C_0((\Gamma,v),M_k)$, for some
  pointed graph $(\Gamma,v)$ and for some $k \in \N$, a \sh{} $\varphi
  \colon B \to A_\infty$, and elements $b_1,b_2,
  \dots, b_n$ in $B$, such that $\|a_j-\varphi(b_j)\| < \ep$ for all $j$.
\end{enumerate}
If, in addition, $A$ is isomorphic to $A \otimes \cO_2$, then
\emph{(i)--(iii)} are equivalent to \emph{(iv)} below, and moreover, the
\sh{} $\varphi$ in \emph{(ii)} and \emph{(iii)} can be taken to be
injective and so can the connecting maps $\varphi_j$ in \emph{(i)}.
\begin{itemize}
\item[\rm{(iv)}]  For every finite set $a_1, a_2, \dots, a_n$ in $A$
  and for every
  $\ep >0$ there exist a sub-\Cs{} $B$ of $A_\infty$ and elements $b_1,b_2,
  \dots, b_n$ in $B$ such that $\|a_j-b_j\| < \ep$ for all $j$ and $B$
  is isomorphic to $C_0(X,M_k)$ for some locally compact
  Hausdorff space $X$ and for some natural number $k$.
\end{itemize}
\end{theorem}

\begin{proof} The implications (i) $\Rightarrow$ (ii) $\Rightarrow$
  (iii) $\Rightarrow$ (iv) are trivial.

Assume that (iv) holds and that $A \cong A \otimes
  \cO_2$. Then $A$ cannot contain a non-zero projection, and as every
  non-zero projection in $A_\infty$ lifts to a non-zero projection in
  $\ell^\infty(A)$, there are no non-zero projections in $A_\infty$.
  Indeed, any
  non-zero projection in $A$ would be infinite (because $A$ is purely
  infinite), and, by (iv), any infinite projection in $A$ is 
  equivalent to an infinite projection
  in a sub-\Cs{} $B$ of $A$ where $B$ is isomorphic to
  $C_0(X,M_k)$ (or to $C(X,M_k)$ if $X$ is compact). But $C_0(X,M_k)$
  (and $C(X,M_k)$) contain no infinite projections.  
 
Let $a_1, \dots, a_n$ in $A$ and $\ep >0$ be given. Find $B \subseteq
A_\infty$ and $c_1, \dots, c_n \in B$ with $B \cong C_0(X,M_k)$ and
$\|c_j-a_j\| < \ep/2$. Use Proposition~\ref{prop:perturbation} to find
a \Cs{} $B_1 \subseteq A_\infty$ and $b_1, \dots, b_n \in B_1$ with $B_1 \cong
C_0((\Gamma,v),M_k)$, or with $B_1 \cong C(\Gamma,M_k)$, and
$\|c_j-b_j\| < \ep/2$, for some closed subset $\Gamma$ 
of $\Gamma_{\ep/2,m}$, some $m$, and some $v \in \Gamma$. We have observed
that $A_\infty$ has no non-zero projections, so $B_1$ has no non-zero
projections. This entails that $\Gamma$ must be connected and it
excludes the unital possibility: $B_1 \cong C(\Gamma,M_k)$. Any
connected closed subset of $\Gamma_{\ep/2,m}$ is (homeomorphic to) a
graph. We have therefore proved that (iv) implies (iii) and that
$\varphi$ in (iii) can be taken to be injective (actually the
inclusion mapping) when $A$ absorbs $\cO_2$. 

Still assuming that $A$ is $\cO_2$-absorbing, one can in (i) and (ii)
arrange that $\varphi$ (in (ii)) and  
the connecting maps $\varphi_j$ (in (i)) are injective. Indeed,
replacing $B$ and $A_j$ by $B/\Ker(\varphi)$ and
$A_j/\Ker(\varphi_{\infty,j})$, respectively, where 
$\varphi_{\infty,j} \colon A_j \to A$ is the inductive limit map, it follows
that $B$ and $A_j$ remain of the form $C_0((\Gamma,v),M_k)$,
with $\Gamma$ a connected graph, because $A$ (and hence $B$ and $A_j$,
respectively) have no non-zero projections.

(iii) $\Rightarrow$ (ii) $\Rightarrow$ (i). These implications follow
by the work of Loring on semiprojective \Cs s. A
finitely generated \Cs{} is semiprojective if and only if it has a
finite set of stable generators (Loring, \cite[Theorem
14.1.4]{Lor:lifting}). It
follows from \cite[Theorem~5.1]{Lor:stableII} (applied with $n=1$)
together with \cite[Theorem 14.1.7]{Lor:lifting} that
$C_0(\Gamma,v)$ is semiprojective when $(\Gamma,v)$ is a pointed
graph. By \cite[Theorem 14.2.2]{Lor:lifting}, any
matrix algebra over $C_0(\Gamma,v)$ is semiprojective.

To prove (iii) $\Rightarrow$ (ii), let $a_1, \dots, a_n \in A$ and
$\ep>0$ be given. Take $\varphi \colon B \to A_\infty$ and $b_1,
\dots, b_n \in B$ as in (iii).   
Let $I_m$ be the ideal in $\ell_\infty(A)$
consisting of those sequences $\{a_k\}_{k=1}^\infty$ for which $a_k=0$
for all $k > m$. Then $\{I_m\}$ is an increasing sequence and $c_0(A) =
\overline{\bigcup_{m=1}^\infty I_m}$. By the definition of
semiprojectivity, for some large enough
  $m$ there is a \sh{} $\psi \colon B \to \ell_\infty(A)/I_m =
  \prod_{k=m+1}^\infty A$ with $\nu_m(\psi(b)) = \varphi(b)$, where $\nu_m
  \colon \ell_\infty(A)/I_m 
  \to A_\omega$ is the quotient mapping. Write $\psi =
  (\psi_{m+1}, \psi_{m+2}, \dots)$, where each $\psi_k \colon
  B \to A$ is a \sh. Now, for $j=1,2, \dots, n$, 
$$\ep > \|a_j-\varphi(b_j)\| = \|\nu_m\big(a_j-\psi_{m+1}(b_j),
a_j-\psi_{m+2}(b_j), \dots\big)\| = \underset{k > m}{\limsup_{k \to
    \infty}} \|a_j - \psi_k(b_j)\|.$$
It follows that $\psi_k \colon B \to A$ satisfies the
conditions of (ii) for some large enough $k$.

The implication (ii) $\Rightarrow$ (i) is contained in
\cite[Lemma~15.2.2]{Lor:lifting}.
\end{proof}

\noindent We conclude this section with a sharpened version of
Proposition~\ref{prop:perturbation}, that will not be used
for the proof of our main result, but is needed for an
application in Section 6. It may also have some independent interest.

\begin{proposition} \label{prop:perturbation1}
Let $X$ be any metrizable compact space. Then there is a metrizable
one-dimensional compact space $Y$ and there is an embedding
$\psi \colon C(Y) \to C(X) \otimes \cO_2$ such that $C(X) \otimes 1
\subseteq \psi(C(Y))$.
\end{proposition}

\begin{proof} To make the strategy of the somewhat lengthy and
  technical proof more transparent, the proof is divided into three
  steps and is carried out backwards. 

\emph{Step 1.} It is shown in Steps 2 and 3 that there is a sequence
of graphs, $\Gamma_n$, and \sh s $\varphi_n \colon C(\Gamma_n) \to C(X)
\otimes \cO_2$ such that 
\begin{itemize}
\item[(a)] $\varphi_1(C(\Gamma_1)) \subseteq \varphi_2(C(\Gamma_2))
\subseteq \varphi_3(C(\Gamma_3)) \subseteq \cdots$,
\item[(b)] $C(X) \otimes 1 \subseteq \overline{\bigcup_{n=1}^\infty
  \varphi_n(C(\Gamma_n))}$.
\end{itemize}
The \sh s $\varphi_n$ need not be injective, but
if we pass to the
quotient of $C(\Gamma_n)$ by the kernel of $\varphi_n$ we obtain
injective \sh s $\psi_n \colon C(Y_n) \to C(X) \otimes \cO_2$
for suitable closed subsets $Y_n$ of $\Gamma_n$, and $\image(\psi_n)
= \image(\varphi_n)$. This leads to the
commutative diagram
$$\xymatrix@C-.3pc{C(Y_1) \ar[d]_{\psi_1} \ar[r]^{\lambda_1} & C(Y_1)
  \ar[d]_{\psi_2} \ar[r]^{\lambda_2} & C(Y_1) \ar[d]_{\psi_3}
  \ar[r]^{\lambda_3} & \cdots \ar[r] & C(Y) \ar@{-->}[d]^-\psi \\
C(X) \otimes \cO_2 \ar@{=}[r] & C(X) \otimes \cO_2 \ar@{=}[r] & C(X)
\otimes \cO_2 \ar@{=}[r] &  \cdots \ar@{=}[r] & C(X) \otimes \cO_2 }
$$
where $\lambda_n = \psi_{n+1}^{-1} \circ \psi_n$ (we are here using
that $\image(\psi_n) \subseteq \image(\psi_{n+1})$). 
The inductive limit of the top row of the diagram above is (isomorphic
to) the Abelian \Cs{} $C(Y)$, where $Y$ is the inverse limit of the
sequence of spaces 
$\{Y_n\}$ with connecting maps $\widehat{\lambda}_n \colon Y_{n+1} \to
Y_n$. In particular, $Y$ is a compact metrizable one-dimensional space
(recall that each of the spaces $Y_n$ is one-dimensional being a
closed subset of the one-dimensional space $\Gamma_n$). The
image of $\psi$ contains $C(X) \otimes 1$ by (b).

\emph{Step 2.} Choose a dense sequence $\{f_1,f_2,f_3, \dots\}$ in
  $C(X)$. The \sh s $\varphi_n$ from Step~1 are constructed from
sequences $\{\varphi^{(k)}_n\}_{k=n}^\infty$ of \sh s
$\varphi^{(k)}_n \colon C(\Gamma_n) \to C(X) \otimes \cO_2$, together
with a collection of finite sets $G_n$ and $F^{(k)}_n$, for $k \ge n$,
all to be constructed in Step~3, with the following properties,
\begin{enumerate}
\item $G_n$ is a set of stable generators, with respect to some
  relations, for the semiprojective \Cs{} $C(\Gamma_n)$,
\item $G_n \subseteq F_n^{(n)} \subseteq F_{n}^{(n+1)} \subseteq
  \cdots \subseteq C(\Gamma_n)$ for all $n$,
\item $\varphi_n^{(k)}(G_n) \subseteq
  \varphi_{n+1}^{(k)}(F_{n+1}^{(k)})$ when $n+1 \le k$,
\item $\|\varphi_{n}^{(k)}(f) - \varphi_{n}^{(k+1)}(f)\| \le
  2^{-k}$ for all $n \le k$ and $f \in F_n^{(k)}$,
\item $\dist(f_j \otimes 1, \varphi_{n}^{(n)}(F_n^{(n)})) \le 1/n$ for all
  $n$ and for $j=1,2, \dots, n$.
\end{enumerate}
It follows from (ii) and (iv) that
\begin{equation} \label{eq:corp1}
\varphi_n(g) \eqdef \underset{k \ge n}{\lim_{k \to \infty}}
\varphi_{n}^{(k)}(g), \qquad g \in G_n,
\end{equation}
exists. The set $\varphi(G_n)$ satisfies the relations satisfied by
$G_n$, whence $\varphi_n$ extends to a \sh{} $\varphi_n
\colon C(\Gamma_n) \to C(X) \otimes \cO_2$ by (i). We show that (a) and (b)
in Step~1 hold, and to do so we first verify the two conditions below:
\begin{itemize}
\item[(c)] $\varphi_n(f) = \lim_{k \to \infty}
\varphi_{n}^{(k)}(f)$ for $f \in C(\Gamma_n)$,
\item[(d)] $\|\varphi_n(f) -
\varphi^{(k)}_n(f)\| \le 2^{-k+1}$ for $f \in F^{(k)}_n$,
\end{itemize}
The set of $f$ for which (c) holds is a \Cs, and this \Cs{} contains
the generating set $G_n$, hence all of $C(\Gamma_n)$.
Property (d) follows easily from (c) and (iv). To see that (a) holds, let $n$
and $g \in G_n$ be
given. For any $k > n$ there is by (iii) an
$h_k \in F_{n+1}^{(k)}$ such that $\varphi_n^{(k)}(g) =
\varphi_{n+1}^{(k)}(h_k)$; whence
$$\|\varphi_n(g) - \varphi_{n+1}(h_k)\| \le \|\varphi_n(g)-
\varphi_n^{(k)}(g)\| + \|\varphi^{(k)}_{n+1}(h_k)-\varphi_{n+1}(h_k)\| \le
2^{-k+2}.$$
This proves that $\varphi_n(G_n)$, and hence
$\varphi_n(C(\Gamma_n))$, are contained in
$\varphi_{n+1}(C(\Gamma_{n+1}))$. To prove (b), take natural numbers
$n \ge j$ and use (v) to find $h_n \in F_n^{(n)}$ such that
$$\|f_j \otimes 1 - \varphi_{n}(h_n) \| \le \|f_j \otimes 1 -
\varphi^{(n)}_n(h_n) \|  +  \|\varphi^{(n)}_n(h_n) -
\varphi_n(h_n)\| \le 1/n + 2^{-n+1}.
$$
(b) clearly follows from this.

\emph{Step 3.}
We construct the graphs
$\Gamma_n$, the \sh s $\varphi_n^{(k)}
\colon C(\Gamma_n) \to C(X) \otimes \cO_2$ and
the sets $F_n^{(k)} \subseteq C(\Gamma_n)$ so that (i)--(v) in Step~2
are satisfied; and this is done by induction on $k$. Since $\cO_2 \cong
\bigotimes_{r=1}^\infty \cO_2$ we can find a sequence of \Cs s
$\C = D_0 \subseteq D_1 \subseteq D_2 \subseteq \cdots \subseteq \cO_2$ with
$D_n \cap D_{n-1}' \cong \cO_2$. The \sh{} $\varphi_n^{(k)}$ will be
constructed so that its image is contained in
$C(X) \otimes D_k$ for all $n \le k$.

For $k=n=1$ use Proposition~\ref{prop:perturbation} to
find an injective \sh{} $\varphi \colon B \to C(X) \otimes D_1$ such that 
$B = C(\Gamma')$ or $B
= C_0(\Gamma',v)$ for some closed subset $\Gamma'$ of
$\Gamma_{1/2,m}$ (for some $m$), and
$\|1 - \varphi(b_1)\| < 1/2$ and $\|f_1 \otimes 1 -\varphi(b_2)\| \le
1$ for some $b_1,b_2 \in B$. We infer from $\|1 - \varphi(b_1)\| <
1/2$ that $B$ is unital and that $\varphi$ is unit-preserving. Hence
$B = C(\Gamma')$. Now, $\Gamma'$ is contained in a compact connected
subset $\Gamma_1$ of $\Gamma_{1/2,m}$, and $\Gamma_1$ is necessarily
(homeomorphic to) a graph. Take $\varphi_1^{(1)} \colon C(\Gamma_1) \to
C(X) \otimes D_1$ to be the composition of the restriction mapping
$C(\Gamma_1) \to C(\Gamma')$ and $\varphi \colon C(\Gamma') \to C(X)
\otimes D_1$. Then $\|f_1 \otimes 1 -\varphi_1^{(1)}(h)\| \le
1$ for some $h \in C(\Gamma_1)$. Let $G_1$ be a set of stable
generators for the semiprojective \Cs{} $C(\Gamma_1)$ with respect to
some relations (cf.\ Loring \cite[Theorem~5.1]{Lor:stableII}), and set
$F_1^{(1)} = G_1 \cup \{h\}$.

Suppose now that $k \ge 1$ and that $\Gamma_n$, $\varphi_n^{(k)}$,
$G_n$, and $F_n^{(k)}$ have been found meeting conditions (i)--(v) 
for all $n = 1, 2, \dots, k$. Then, by (i) and (iii),
$$\varphi_1^{(k)}(C(\Gamma_1)) \subseteq \varphi_2^{(k)}(C(\Gamma_2))
\subseteq \varphi_3^{(k)}(C(\Gamma_3)) \subseteq \cdots \subseteq
\varphi_k^{(k)}(C(\Gamma_k)) \subseteq C(X) \otimes D_k.$$
We proceed to find $\Gamma_{k+1}$ and \sh s
$$\varphi_{k+1}^{(k+1)} \colon C(\Gamma_{k+1}) \to C(X) \otimes
D_{k+1}, \quad \varphi_n^{(k+1)} \colon C(\Gamma_n) \to
\varphi_{n+1}^{(k+1)}(C(\Gamma_{n+1})), \; \,  1 \le n \le k,$$
such that
$$\dist(f_j \otimes 1, \image(\varphi_{k+1}^{(k+1)})) \le
1/(k+1), \qquad \|\varphi_n^{(k+1)}(f)-\varphi_n^{(k)}(f)\| \le
2^{-k},$$
for $1 \le j \le k+1$, $1 \le n \le k$, and for all $f \in
F_n^{(k)}$. Once this is done, it is straightforward to find the sets
$F_n^{(k+1)}$, $n=1, \dots, k+1$, such that (ii)--(v) are satisfied.

As the generating set $G_n$ for
$C(\Gamma_n)$ is
stable, $n=1, \dots, k$, one can for each $\ep>0$ and for each finite
subset $F$ of
$C(\Gamma_n)$ find $\delta = \delta(\ep,F,\Gamma_n, G_n) >0$ such that
whenever $B \subseteq A$ are
unital \Cs s and $\rho \colon C(\Gamma_n) \to A$ is a \sh{} with
$\dist(\rho(g),B) \le \delta$ for all $g \in G_n$, then there is a
\sh{} $\sigma \colon C(\Gamma_n) \to B$ with $\|\rho(f)-\sigma(f) \| \le
\ep$ for all $f \in F$.

Set $\delta_1 = 2^{-k}$ and set
\begin{eqnarray*}
\delta_{n+1} & = & \delta(\delta_{n},F_n^{(k)}, \Gamma_n,G_n) \wedge 2^{-k},
\quad n = 1,2, \dots, k.
\end{eqnarray*}
Arguing as in the case $k=n=1$ we use
Proposition~\ref{prop:perturbation} to find a graph $\Gamma_{k+1}$ and
a \sh{} $\varphi_{k+1}^{(k+1)} \colon C(\Gamma_{k+1}) \to C(X) \otimes
D_{k+1}$ such that
$$\dist(f_j \otimes 1, \image(\varphi_{k+1}^{(k+1)})) \le 1/(k+1),
\qquad \dist(\varphi_k^{(k)}(g), \image(\varphi_{k+1}^{(k+1)}))
\le \delta_{k+1},$$
for $1 \le j \le k+1$ and for all $g \in G_k$. Take a finite set
$G_{k+1}$ of stable generators for the semiprojective \Cs{}
$C(\Gamma_{k+1})$. By
the choice of $\delta_{k+1}$ we can then find a \sh{}
$\varphi_k^{(k+1)} \colon C(\Gamma_k) \to
\image(\varphi_{k+1}^{(k+1)}) \subseteq C(X) \otimes \cO_2$ such that 
$\|\varphi_k^{(k+1)}(f) -  \varphi_k^{(k)}(f)\| \le \delta_k \; (\le
2^{-k})$ for all $f$ in $F_k^{(k)}$. In particular, 
$\dist(\varphi_{k-1}^{(k)}(g), \image(\varphi_{k}^{(k+1)}))
\le \delta_{k}$ for all $g \in G_{k-1}$. We can in this way continue the
construction of $\varphi_{k-1}^{(k+1)},\varphi_{k-2}^{(k+1)}, \dots,
\varphi_{1}^{(k+1)}$ with the desired properties. 
\end{proof}

\noindent There is a curious topological corollary to
Proposition~\ref{prop:perturbation1} that each compact metrizable
space is the continuous-open image of a one-dimensional
space. This relates to the classical Peano Curve: a continuous
surjection from the 
interval onto the square (or onto more arbitrary compact spaces). The
Peano Curve is not open, and there exists no open continuous
surjection from the interval onto the square. 
Indeed, suppose that $f$ were a continuous and open mapping
from the interval onto the square, let $I$ be a closed line segment in the
square, and let $C = f^{-1}(I)$. The interior of $C$ is empty because
the interior of $I$ is empty, so $C$ is compact and totally
disconnected. The restriction mapping $f \colon C \to I$ is continuous,
open, and maps compact sets to compact sets; hence it maps clopen sets
into clopen sets. Thus $f(C_0)=I$ 
for every non-empty clopen subset of $C_0$ of $C$. As $C$ has clopen
subsets of arbitrary small diameter, this contradicts the continuity
of $f$.  

We thank Etienne Blanchard for suggesting the continuous bundle proof
below.

\begin{corollary} \label{cor:Peano}
For each compact metrizable space $X$ there is a compact metrizable
one-dimensional space $Y$ and a continuous, open, surjection $Y \to X$.
\end{corollary}

\begin{proof}
Use Proposition~\ref{prop:perturbation1} to find a compact metrizable
one-dimensional space $Y$ and an embedding $\varphi \colon C(Y) \to
C(X) \otimes \cO_2$ such that $C(X) \otimes 1 \subseteq
\varphi(C(Y))$. Then $f \mapsto \varphi^{-1}(f \otimes 1)$ defines an
injection from $C(X)$ into $C(Y)$ which can be realized by a
continuous surjective function $\lambda \colon Y \to X$; i.e., $f
\circ \lambda = \varphi^{-1}(f \otimes 1)$, or $\varphi(f \circ
\lambda) = f \otimes 1$, for all $f$ in $C(X)$.

We proceed to show that $\lambda$ is open. The argument uses
implicitly that $C(Y)$ is a continuous field over $X$. Let $\pi_x
\colon C(X) \otimes \cO_2 \to \cO_2$ be the point evaluation at
$x$, and put $\nu_x(g) = \pi_x(\varphi(g))$ for $g \in C(Y)$. Then $x
\mapsto \|\nu_x(g)\|$ is continuous for every $g$ in $C(Y)$. We show that
$$\|\nu_x(g)\| = \sup\{|g(y)| : y \in \lambda^{-1}(x) \}, \qquad x \in
X, \; g \in C(Y).$$
Fix $x$ in $X$ and choose a decreasing sequence
$\{f_n\}_{n=1}^\infty$ of positive contractions in $C(X)$ such that
$f_n(x) = 1$ for all $n$ and $f_n(x') = 0$ eventually when $x' \ne x$. Then
$$\|\nu_x(g)\| = \lim_{n \to \infty}  \|(f_n\otimes 1)\varphi(g)\| =
\lim_{n \to  \infty}  \|(f_n \circ \lambda)g\| =  \sup\{|g(y)| : y \in
\lambda^{-1}(x) \}.$$

Let $U$ be an open subset of $Y$. Because $Y$ is a metrizable space there
is a function $g \in C(Y)$ such that $g(y) \ne 0$ if and only if $y
\in U$. As
$$x \in \lambda(U) \iff \lambda^{-1}(x) \cap U \ne \emptyset \iff
g|_{\lambda^{-1}(x)} \ne 0 \iff \|\nu_x(g)\| >0,$$
it follows that $\lambda(U) = \{x \in X : \|\nu_x(g)\| > 0\}$, and
this proves that $\lambda(U)$ is open.
\end{proof}

\section{A Hahn--Banach separation theorem for \cp{} maps}
\label{sec:cp}

\noindent  This section contains a Hahn--Banach type separation
theorem for \cp{} maps, a result implicitly contained in
\cite{KirRor:pi2}, and some applications.

\begin{definition} \label{def:convex}
Let $A$ and $B$ be \Cs s. Denote by  $CP(A,B)$ the cone of all
completely positive maps from $A$ to $B$. A subset $\sK$ of $CP(A,B)$
is called \emph{operator convex} if it satisfies:
\begin{enumerate}
\item $\sK$ is a cone.
\item If $V \in \sK$ and $b \in B$, then
  the map $a \mapsto b^*V(a)b$ belongs to $\sK$.
\item If $V \in \sK$, $c_1, \dots, c_n \in A$, and
  $b_1, \dots, b_n \in B$, then the map given by
$$a \mapsto \sum_{i=1}^n \sum_{j=1}^n b_i^*V(c_i^*ac_j)b_j,$$
belongs to $\sK$.
\end{enumerate}
\end{definition}

\noindent The \cp{} map displayed in (iii) is equal to
the map $a \mapsto  b^*(V \otimes \Id_{M_n})(c^*ac)b$, where $b$ is
the row matrix $(b_1, \dots, b_n)$ and $c$ is the column matrix
$(c_1,\dots, c_n)^t$.

\begin{proposition}[cf.\ Lemma 7.18 of \cite{KirRor:pi2}] \label{prop:HB}
Let $A$ be a separable, nuclear \Cs{} and let $B$ be any \Cs. Suppose that
$\sK$ is an operator convex subset of $CP(A,B)$, and let $T \in
CP(A,B)$. Then $T$
belongs to the point-norm closure of $\sK$ if and only if $T(a)$
belongs to the closed two-sided ideal in $B$ generated by
$\{V(a) \mid V \in \sK \}$ for every $a$ in $A$.
\end{proposition}

\begin{proof} The ``only if'' part is trivial. The proof of the ``if''
  part is almost verbatim identical with the proof of
  \cite[Lemma~7.18]{KirRor:pi2}, but since there are some changes (in
particular in  the last paragraph of the proof and in the notation),
we repeat below most of the argument for the convenience of the reader.

For each finite set $a_1, \dots, a_n$ in $A$ and for each $\ep>0$ we
must find $V$ in $\sK$ such that $\|T(a_j)-V(a_j)\|
\le \ep$ for all $j$. Equivalently, we must show that
\begin{equation} \label{eq:HB1}
(t(a_1), \dots, T(a_n)) \in \overline{\{\big(V(a_1), \dots,
  V(a_n)\big) \mid V \in \sK \}}.
\end{equation}
By a Hahn--Banach separation argument, and because the set on the
right-hand side of \eqref{eq:HB1} is a cone (by assumption (i)), and
hence convex, \eqref{eq:HB1} is implied by the following: For each
set $f_1, \dots, f_n$ in $B^*$ and for
each $\ep>0$ there is $V$ in $\sK$ such $|f_j(T(a_j)) - f_j(V(a_j))|
\le \ep$ for $j=1, \dots, n$.

Choose a cyclic representation $\pi
\colon B \to B(H)$, a cyclic
vector $\xi$ in $H$, and elements $c_1, \dots, c_n$ in $\pi(B)' \cap
B(H)$ such that $f_j(a) = \langle \pi(a)\xi, c_j^*\xi \rangle$ for all
$a$ in $A$ and for all $j$; cf.\
\cite[Lemma~7.17~(i)]{KirRor:pi2}.  Let $C$ be the
\sCs{} of $\pi(B)' \cap B(H)$ generated by $c_1, \dots,
c_n$. Keeping $\pi$ and $\xi$ fixed, associate to each \cp{} map $V
\colon A \to B$,  the positive functional $\varphi_V$ on $A \otimes C$
given by
\begin{equation} \label{eq:nuclear}
\varphi_V(a \otimes c) = \big\langle \pi(V(a)) \xi, c^* \xi \big\rangle,
\qquad a \in A, \; \; c \in C.
\end{equation}
(A priori, $\varphi_V$ defines a functional on the maximal tensor
product $A \otimes_{\mathrm{max}} C$, but the maximal and the minimal
tensor products on $A \odot C$ coincide because $A$ is nuclear.) Let
$\sC$ be the weak-$^*$ closure of the cone $\{\varphi_V : V \in
\sK \}$. Observe that $f_j(V(a_j)) = \langle \pi(V(a_j)) \xi, c_j^*
\xi \rangle = \varphi_V(a_j \otimes c_j)$ for every \cp{} map $V
\colon A \to B$. Hence it will suffice to show that $\varphi_T$
belongs to $\sC$.

It follows as in the proof of  \cite[Lemma~7.18]{KirRor:pi2} (using
assumption (iii)!) that $d^*\rho d$ belongs to $\sC$ for all
$\rho$ in $\sC$ and for all $d$ in $A \otimes C$.

As in \cite[Lemma~7.17~(ii)]{KirRor:pi2}, let $J$ be the
closed two-sided ideal
in $A \otimes C$ consisting of those elements $z$ for which
$\varphi_V(z^*z)=0$ for all $V$ in $\sK$. By
\cite[Lemma~7.17~(ii)]{KirRor:pi2} one can conclude that
$\varphi_T$ belongs to $\sC$ if we know that $\varphi_T(z^*z)=0$ for
all $z$ in $J$. Using again the assumption that $A$ is nuclear,
we can use a theorem of Blackadar, \cite[Theorem 3.3]{Bla:tensor}
(see also \cite[Proposition 2.13]{Kir:fields}) that $J$ is the closed
linear span
of the set of elementary tensors $x \otimes y$ in $J$. The left
kernel $L$ of $\varphi_T$, consisting of all $z$ in $A \otimes C$
such that $\varphi_T(z^*z)=0$, is a closed linear subspace of $A
\otimes C$, and so it
suffices to show that $\varphi_T(x^*x \otimes y^*y)=0$ whenever $x
\in A$ and $y \in C$ are such that $x \otimes y$ belongs to $J$

Fix $x \otimes y$ in $J$ as above. By assumption, $T(x^*x)$ belongs to
the closed two-sided ideal in $B$
generated by $\{V(x^*x) \mid V \in \cK\}$. For any
$\ep >0$ we can therefore find $V_1,
\dots, V_m$ in $\sK$ and $b_1, \dots, b_m$ in $B$ such that
$$\|T(x^*x) - b_1^*V_1(x^*x)b_1 - \cdots - b_m^*
V_m(x^*x)b_m\| \le \ep.$$
In other words, $\|T(x^*x) - W(x^*x)\| \le \ep$, when $W(a) =
b_1^*V_1(a)b_1 + \cdots + b_m^*V_m(a)b_m$. By assumptions (i) and (ii)
we see that $W$ belongs to $\sK$. Hence,
\begin{eqnarray*}
|\varphi_T(x^*x \otimes y^*y)| & = & |\varphi_T(x^*x \otimes y^*y) -
\varphi_W(x^*x \otimes y^*y)| \\ & = &
\big|\big\langle \pi(T(x^*x) - W(x^*x))\xi, y^*y \xi
\big\rangle \big| \; \le \;  \ep\|\xi\|^2\|y\|^2.
\end{eqnarray*}
As $\ep>0$ was arbitrary, we get $\varphi_T(x^*x \otimes y^*y)=0$
as desired.
\end{proof}

\noindent The nuclearity assumption in
Proposition~\ref{prop:HB} is necessary. Indeed, for any pair
of \Cs s $A$ and $B$, the set
$CP_{\mathrm{nuc}}(A,B)$ of all nuclear \cp{} maps from $A$ to $B$ is
operator convex; $a = \Id_A(a)$
belongs to the closed two-sided ideal generated by $\{V(a) \mid V \in
CP_{\mathrm{nuc}}(A,A)\}$ for every $a \in A$; but  $\Id_A$
belongs to the point-norm closure of $CP_{\mathrm{nuc}}(A,A)$ if and
only if $A$ is nuclear.

To prove the second claim, take a non-zero $a \in A$, take a
positive linear functional $\rho$ on $A$ such that $\rho(a) \ne 0$, and
set $V(x) = \rho(x)a^*a$, $x \in A$. Then $V$ belongs to
$CP_{\mathrm{nuc}}(A,A)$
(because $V$ factors through the complex numbers via the positive
functional $\rho$), and $a$ belongs to the closed two-sided ideal in
$A$ generated by $V(a) = \rho(a)a^*a$.

\begin{corollary} \label{cor:appr_inner}
Let $B$ be any \Cs, let $A$ be a separable nuclear sub-\Cs{} of $B$,
and let $T \colon A \to B$ be a \cp{} map. Then $T$ is
approximately inner (cf.\ Definition~\ref{def:app_inner})
if and only if $T(a) \in \overline{BaB}$ for all $a \in A$.
\end{corollary}

\begin{proof} The set, $\sK$, of all inner
  \cp{} maps from $A$ into $B$ is operator convex. Indeed, it is
  trivial that (i) and
  (ii) hold of Definition~\ref{def:convex} hold. To see that (iii)
  holds, take $V$ in $\sK$, take elements
  $c_1, \dots, c_m \in A$, take elements $b_1, \dots, b_m \in B$, and
  set $W(a) = \sum_{i,j=1}^n b_i^*V(c_i^*ac_j)b_j$. If $V(a) =
  \sum_{k=1}^n d_k^*ad_k$, then $W(a) = \sum_{k=1}^n f_k^*af_k$, when
  $f_k = \sum_{j=1}^n c_jd_kb_j$.

Observe that $T \colon A \to B$ is approximately inner if and only if
$T$ belongs to the point-norm closure of $\sK$. The set $\{V(a) \mid V
\in \sK\}$ is contained in the two-sided ideal in $B$ generated by $a$ (it
actually coincides with the positive part of the algebraic two-sided
ideal generated by $a$, when $a$ is positive), and so the lemma follows from
Proposition~\ref{prop:HB}.
\end{proof}

\noindent We do not know if the conclusion of
Corollary~\ref{cor:appr_inner} actually characterizes nuclear \Cs s,
as phrased more formally below:

\begin{question} \label{q4} Suppose that $A$ is a \Cs{} for which
  every \cp{} map $V \colon A \to A$, that satisfies $V(a) \in
  \overline{AaA}$ for all $a \in A$, is approximately inner. Does it
  follow that $A$ is nuclear?
\end{question}

\noindent We conclude this section with a (new) proof of the
fundamental uniqueness theorem for $\cO_2$-absorbing \Cs s. Two \sh s
$\varphi, \psi \colon A \to B$ are said to induce the same map at the
level of ideals, if $\varphi^{-1}(I) = \psi^{-1}(I)$ for all closed
two-sided ideals $I$ in $B$. This is the case if and only if the closed
two-sided ideals in $B$ generated by
$\varphi(a)$ and $\psi(a)$, respectively, coincide for all $a \in
A$.

\begin{theorem} \label{thm:O_2-uniqueness}
Let $A$ be a separable, nuclear \Cs, let $B$ be an
$\cO_2$-absorbing stable \Cs, and let $\varphi, \psi \colon A \to B$
be \sh s. Then $\varphi$ and $\psi$ induce the
same map at the level of ideals if and only if they are
approximately unitarily equivalent,
i.e., there is a sequence $\{u_n\}_{n=1}^\infty$ of unitary elements
of $\cM(B)$ such that $u_n^*\psi(a)u_n \to \varphi(a)$ for all $a
\in A$.
\end{theorem}

\begin{proof} The ``if'' part is trivial. To prove the ``only if''
  part, let $\sK$ denote the set of all \cp{} maps $A \to B$ of the
  form $V \circ \varphi$, where $V \colon \varphi(A) \to B$ is approximately
  inner. Then $\sK$ is operator convex
  (cf.\ the proof of Corollary~\ref{cor:appr_inner}). By assumption,
  $\psi(a)$ belongs to 
  $\overline{B\varphi(a)B}$, and this ideal is contained in the closed
  two-sided ideal generated by $\{W(a) \mid W \in
  \sK\}$ for all $a \in A$, whence $\psi$ belongs to the point-norm
  closure of $\sK$ by Proposition~\ref{prop:HB}.  It follows
  from \cite[Proposition~8.4]{KirRor:pi2} that all approximately inner
  \cp{} maps $\varphi(A) \to B$ are approximately 1-step
  inner. Accordingly, there is a sequence $\{d_n\}$ in $B$
  such that $d_n^*\varphi(a)d_n \to \psi(a)$ for all $a \in A$.

Let $\{e_n\}_{n=1}^\infty$ be an approximate unit for $A$ consisting
of positive contractions, and put $c_n = \varphi(e_n)d_n$. Then
$c_n^*c_n = d_n^*\varphi(e_n^2)d_n \to \psi(e_n^2)$, whence $\|c_n\|
\to 1$, and $c_n^*\varphi(a)c_n \to \psi(a)$ for all $a \in A$. 
The contractions $b_n = \|c_n\|^{-1}c_n$ satisfy $b_n^*\varphi(a)b_n
\to \psi(a)$ for all $a \in A$. By \cite[Lemma~2.4]{KirRor:pi2} there
is a sequence $\{s_n\}$ of isometries in $\cM(B)$ such that
$s_n^*\varphi(a)s_n \to \psi(a)$ for all $a \in A$. 

A similar argument shows that there is a sequence $\{t_n\}$ of
isometries in
$\cM(B)$ satisfying $t_n^*\psi(a)t_n \to \varphi(a)$ for all $a \in A$.

For all $x,y \in \varphi(A)$,
$$\lim_{n \to \infty} \|s_n^*xys_n - s_n^*xs_ns_n^*ys_n\| = 0, \qquad 
\lim_{n \to \infty} \|t_n^*s_n^*xs_nt_n - x\| = 0,
$$
and so it follows from Lemma~\ref{lm:cantor-bernstein} that there is a
sequence $\{u_n\}$ of unitary elements in $\cM(B)$ such that
$u_n^*\psi(a)u_n \to \varphi(a)$ for all $a \in A$. 
\end{proof} 

\section{The main result: ideal-system preserving zero homotopic
  purely infinite \Cs s}
\label{sec:homotopy-trivial}

\noindent In this section we prove our main result. As we shall spend a great
deal of effort to keep track of ideals, it is convenient to set up
some relevant notation. If $F$ is a subset of a \Cs{} $A$, then
$\overline{AFA}$ and $I_A(F)$ will both denote the closed two-sided
ideal in $A$ generated by $F$. In the case where $F = \{a\}$ is a
singleton, we write $\overline{AaA}$ and $I_A(a)$, respectively, instead
of $\overline{AFA}$ and $I_A(F)$.

\begin{definition} \label{def:homotopy-trivial}
A \Cs{} $A$ is said to be \emph{zero-homotopic in an ideal-system preserving
  way} if there is a pointwise-norm continuous path $\{\rho_t\}_{t
  \in [0,1]}$ of $^*$-endomorphisms on $A$ such that $\rho_0 = 0$,
$\rho_1 = \Id_A$, and $\rho_t(J) \subseteq J$ for all closed
two-sided ideals $J$ in $A$.
\end{definition}

\noindent A word of warning: The cone $CA = C_0((0,1],A)$ over a \Cs{}
$A$ is usually not zero-homotopic in an ideal-system preserving
  way (although it is zero-homotopic). The standard homotopy path
  $\varphi_t \colon CA \to CA$ from the identity to zero is not ideal
  preserving. In Section~\ref{sec:applications} we shall give
  examples of \Cs s that are  zero-homotopic in an ideal-system preserving
  way.

We emphasize below an important consequence of the classification of
separable, nuclear, $\cO_\infty$-absorbing, stable \Cs s by a an ideal
respecting version of Kasparov's $KK$-theory.

\begin{theorem}[Kirchberg \cite{Kir:fields}] \label{thm:Kir}
Let $A$ be a separable, stable, $\cO_\infty$-absorbing, nuclear \Cs{} that is
zero-homotopic in an ideal-system preserving way. Then $A \cong A
\otimes \cO_2$.
\end{theorem}

\begin{lemma} \label{lm:C}
Let $A$ be a separable stable \Cs. Then $A$ contains an Abelian
sub-\Cs{} $C$ that separates ideals of $A$.
\end{lemma}

\noindent That $C$ separates ideals of $A$ means that $I \cap C = J
\cap C \Rightarrow I = J$ for any two
closed two-sided ideals $I,J$ of $A$.

\begin{proof} Because $A$ is stable and separable we can write $A =
  A_0 \otimes \cK$ for some separable \Cs{} $A_0$. Let $\{a_1,a_2,
    \dots\}$ be a dense subset of the positive cone in $A_0$. Let
    $\{p_{n,m}\}_{n,m=1}^\infty$ be pairwise orthogonal 1-dimensional
    projections in $\cK$, and put
$$C = C^*\big((a_n-1/m)_+ \otimes p_{n,m} \mid n,m \in \N\big).$$
Clearly, $C$ is an Abelian sub-\Cs{} of $A$. We must check that $C$
separates ideals, so let $I$ and $J$ be closed two-sided ideals in
$A=A_0 \otimes \cK$ with $I \cap C
= J \cap C$.  Note that $I = I_0 \otimes \cK$ and $J = J_0 \otimes
\cK$ for some ideals $I_0$ and $J_0$ in $A_0$. Let $a$ be a positive
element in $I_0$ and let $\ep>0$. Find $n,m \in \N$ such that
$\|a-a_n\| < 1/m < \ep$. Then $(a_n-1/m)_+ = x^*ax$ for some $x$ in
$A_0$ (see \cite[Lemma~2.2]{KirRor:pi2}), whence
$$(a_n-1/m)_+ \otimes p_{n,m} \in I_A(a \otimes p_{n,m}) \cap C
\subseteq I \cap C = J \cap C,$$
and so $(a_n-1/m)_+$ belongs to $J_0$. As
$$\|a-(a_n-1/m)_+\| \le \|a-a_n\| + \|a_n - (a_n-1/m)_+\| < 1/m + 1/m
\le 2 \ep,$$
and because $\ep >0$ was arbitrary, it follows that $a$ belongs to
$J_0$. This proves $I_0 \subseteq J_0$, and a symmetric argument
yields the reverse inclusion.
\end{proof}

\noindent
The reader is referred to \cite[Definition~5.1]{KirRor:pi2}
for a formal definition being ``\emph{strongly purely infinite}''. Every
$\cO_\infty$-absorbing \Cs{} is strongly purely infinite, and it is
proved in \cite{KirRor:pi2} that every separable, nuclear, strongly
purely infinite \Cs{} that is either unital or stable is
$\cO_\infty$-absorbing.

\begin{proposition}[Proposition 7.13 of \cite{KirRor:pi2}]
\label{prop:extension}
Suppose that $A$ is a separable strongly purely infinite \Cs. Let
$C$ be an Abelian sub-\Cs{} of $A$ and let $\omega$ be a free filter
on $\N$. Then there is an Abelian \Cs{} $C_1$ with $C \subseteq C_1
\subseteq A_\omega$ and a 1-step inner \cpc{} $V \colon A \to C_1$ with
$V(c) =c$ for all $c$ in $C$.
\end{proposition}

\noindent From now on, let $V \colon A \to C_1 \subseteq A_\omega$ be
the \cpc{} as in Proposition~\ref{prop:extension} wrt.\
the Abelian sub-\Cs{} $C$ found in Lemma~\ref{lm:C}.

\vspace{.3cm}
\noindent {\bf{The Kasparov--Stinespring dilation.}}
Kasparov's version of the Stinespring dilation,
\cite[Theorem~3]{Kas:cp}, shows that every \cpc{} $T \colon A \to B$ has a
dilation to a \sh{} $\psi_T \colon A \to \cM(B \otimes \cK)$ which
satisfies
\begin{equation} \label{eq:Kasparov}
(1 \otimes e_{11})\psi_T(a)(1 \otimes e_{11}) = T(a) \otimes e_{11},
\qquad a \in A,
\end{equation}
where $e_{11}$ is a 1-dimensional projection in $\cK$. We shall apply
this to the \cpc{} $V
\colon A \to C_1$ considered above.

Because $V(c)=c$ for all $c \in C \subseteq C_1$ the two mappings $c
\mapsto \psi_V(c)$ and $c \mapsto c \otimes 1$ from $C$ into
$\cM(C_1 \otimes \cK)$ are approximately unitarily equivalent. In
particular, $I_{\cM(C_1\otimes\cK)}(\psi_V(c)) =
I_{\cM(C_1\otimes\cK)}(c \otimes 1)$ for all $c \in C$.
We note for later use that $\psi_C$ is ideal preserving as follows:

\begin{lemma} \label{lm:psi_V}
For each closed two-sided ideal $J$ in $A$, set $J_1 = I_{C_1}(J \cap
C)$, and denote by  $\cM(C_1 \otimes
\cK,J_1 \otimes \cK)$ the closed two-sided ideal in $\cM(C_1
\otimes \cK)$ consisting of those multipliers $x \in \cM(C_1 \otimes
\cK)$ for which $xy, yx \in J_1 \otimes \cK$ for all $y \in C_1
\otimes \cK$. Then
\begin{equation} \label{eq:psi_V}
\psi_V(J)  \subseteq \cM(C_1 \otimes \cK,J_1 \otimes \cK).
\end{equation}
\end{lemma}

\begin{proof}
Take first $c \in J \cap C$. Then
$$I_{\cM(C_1\otimes\cK)}(\psi_V(c)) \cap C_1 \otimes \cK \; = \;
I_{\cM(C_1\otimes\cK)}(c \otimes 1)  \cap C_1 \otimes \cK \; \subseteq
\; J_1 \otimes \cK,$$
whence $\psi_V(c)$
belongs to $\cM(C_1 \otimes \cK,J_1 \otimes \cK)$. It follows that
$\psi_V(a) \in \cM(C_1 \otimes \cK,J_1 \otimes \cK)$ for all $a \in
I_A(J \cap C)$. But $I_A(J \cap C) = J$ for every closed two-sided
ideal $J$ in $A$ because $C$ separates ideals in $A$.
\end{proof}

\vspace{.3cm}
\noindent We shall from now on consider a \Cs{} $A$ that has the
following properties:
\begin{equation} \tag{$\dagger$}
\text{\emph{separable, strongly purely infinite, and zero-homotopic in
  an ideal preserving way.}}
\end{equation}
Property ($\dagger$) passes to quotients. In particular, $A$ has
no non-zero projections, and no quotient of $A$ contains a non-zero
projection. As $A$ in particular is purely infinite and has no unital
quotient, it follows from \cite[Theorem~4.24]{KirRor:pi} that $A$ is stable.
If $A$ besides having property ($\dagger$) is nuclear, then $A \cong A \otimes
  \cO_\infty$ (by \cite[Theorem~8.6]{KirRor:pi2}). Furthermore, when
  $A$ is nuclear, Theorem~\ref{thm:Kir} implies that $A \cong A
  \otimes \cO_2$.

Fix a path $\{\rho_t\}_{t \in [0,1]}$ of $^*$-endomorphism
that implements the zero-homotopy and define the \sh{}
$$\rho \colon A \to C_0( ] 0,1]
, A) \quad \text{by} \; \; \rho(a)(t) =
\rho_t(a), \; a \in A.$$
Note that $\rho(J) \subseteq C_0( ] 0,1], J)$ for all closed two-sided
ideals $J$ of $A$ by assumption on the path $\rho_t$.

In the following lemmas we show that if $A$ is nuclear (besides
having property ($\dagger$)), then it
satisfies condition (iv) of Theorem~\ref{thm:ah_0} wrt.\ any free
filter $\omega$ on $\N$. (We can take $\omega = \omega_\infty$.)  The
idea of the proof is outlined in the diagram:
\begin{equation} \label{eq:diagram}
\begin{split}
\xymatrix{A \ar[r]^-\rho \ar[d]_{\psi_V}  & C_0( ] 0,1], A)
  \ar@{-->}[d]^-{W_n} & & \\ \cM(C_1 \otimes \cK) & C_1 \otimes \cK
  \ar@{_{(}->}[l] \ar@{-->}[r]_-{\iota} & A_\omega \ar@{-->}[r]_-{T_n}
  & A_\omega
}
\end{split}
\end{equation}
The maps $W_n$, $\iota$, and $T_n$ (that are to be
constructed) are respectively an almost multiplicative \cpc{}, a
\sh{}, and a \cpc{}, chosen such that the sequence $\{T_n \circ \iota \circ
W_n \circ \rho\}$
converges pointwise to the inclusion mapping $A \to A_\omega$. Using
Remark~\ref{rem:cantor-bernstein} (among other results) one finds a
sequence $\{u_n\}$ of
unitaries in $\cM(A)_\omega$
such that $\|u_n  (\iota \circ W_n \circ
\rho)(a)u_n^*-a\| \to 0$ for all $a \in A$. As $C_1 \otimes \cK$ is
the completion of $\bigcup_{k=1}^\infty C_1 \otimes M_k$, we can
choose the sub-\Cs{} $B$ of $A_\omega$ as in
Theorem~\ref{thm:ah_0}~(iv) to be $u_n(\iota(C_1
\otimes M_{k}))u_n^*$  for some large enough $n$ and $k$.

We begin by constructing the mapping $W_n$. For each
natural number $n$ choose a continuous function $g_n \colon [0,1]^2
\to \R^+$ such that
\begin{equation} \label{eq:g}
\int_0^1 g_n(s,t)^2 dt = 1 \; \; \text{for all} \; \; s \in [0,1],
\quad \text{and} \quad g_n(s,t) = 0
\; \; \text{when} \; \; |s-t| \ge 1/n.
\end{equation}
Choose a quasi-central approximate unit $\{e_n\}_{n=1}^\infty$ for $C_1 \otimes
\cK$ consisting of
positive contractions that satisfy $e_n e_{n+1} = e_n = e_{n+1}e_n$
and $e_n(1 \otimes e_{11}) = (1 \otimes e_{11})e_n$ for all $n$, and
$$\lim_{n \to \infty} \sup_{t \in [0,1]} \|g_n(e_n,t)\psi_V(a) -
\psi_V(a)g_n(e_n,t)\| = 0, \qquad  a \in A.$$
For each $f \in C_0(] 0,1], A)$ and for each natural number $n$ put.
\begin{equation*} 
\widetilde{W}_n(f) = \int_0^1 g_n(e_n,t) \psi_V(f(t)) g_n(e_n,t) dt, \quad
W_n(f) = \int_{1/n}^1 g_n(e_n,t) \psi_V(f(t)) g_n(e_n,t) dt,
\end{equation*}

We record below some properties of the maps $W_n$.

\begin{lemma} \label{lm:W_n}
Let $f, g$ be in $C_0(]0,1],A)$ and let $m < n$ be natural numbers.
\begin{enumerate}
\item ${\displaystyle{\lim_{n \to \infty}}} \|\widetilde{W}_n(f)-W_n(f)\| = 0$.
\item $W_n$ is a \cpc{} from $C_0(]0,1],A)$ into $C_1 \otimes \cK$.
\item $W_n(f)$ belongs to the ideal in $C_1 \otimes \cK$ generated by
  $\{e_n \psi_V(f(t))e_n \mid t \in [0,1]\}$.
\item ${\displaystyle{\lim_{n \to \infty}}} \|W_n(fg) -
  W_n(f) W_n(g)\| = 0$.
\item $\|e_m W_n(f)e_m - e_m \psi_V(f(1)) e_m \| \le \sup\{\|f(t)-f(1)\|
  : t \in [1-1/n,1] \}$.
\end{enumerate}
\end{lemma}

\begin{proof} (i). Since $f(t) \to 0$ as $t \to 0$ we get
\begin{eqnarray*}
\|\widetilde{W}_n(f)-W(f)\| & = & \big\|\int_0^{1/n} g_n(e_n,t) \psi_V(f(t))
g_n(e_n,t) dt \big\|  \\ & \le & \sup_{t \in [0,1/n]} \|f(t)\| \cdot
\big\|\int_0^1
g_n(e_n,t)^2 dt \big\| = \sup_{t \in [0,1/n]} \|f(t)\| \to 0.
\end{eqnarray*}

(ii) and (iii). Approximating the integral that defines $W_n$ by a
Riemann sum we see that $W_n$ is a \cp{} map.
Note that $g(0,t) = 0$ for $t \ge 1/n$. It follows that
$g(e_n,t)$ belongs to $C^*(e_n)$ for $t \ge 1/n$, whence
$g(e_n,t)\psi_V(f(t))g(e_n,t)$ belongs 
to the ideal in $C_1 \otimes \cK$ generated by $e_n\psi_V(f(t))e_n$.

(iv). By (i) it will suffice to show that $\|\widetilde{W}_n(fg) -
  \widetilde{W}_n(f) \widetilde{W}_n(g)\| \to 0$.
Note that $g_n(x,t)g_n(x,s) = 0$ if $x,t,s \in [0,1]$  and
$|t-s| \ge 2/n$. It follows that $g_n(e_n,t)g_n(e_n,s)=0$ when $|t-s|
\ge 2/n$. Hence
\begin{eqnarray*}
\widetilde{W}_n(f) \widetilde{W}_n(g) &=& \int_0^1 \int_0^1
g_n(e_n,t)\psi_V(f(t))g_n(e_n,t)g_n(e_n,s)\psi_V(g(s))g_n(e_n,s) \, ds
\, dt
\\ &=& \int_0^1\int_{(t-2/n) \vee 0}^{(t+2/n)\wedge 1}
g_n(e_n,t)\psi_V(f(t))g_n(e_n,t)g_n(e_n,s)\psi_V(g(s))g_n(e_n,s) \, ds
\, dt.
\end{eqnarray*}
By the choice of the quasi-central approximate unit $\{e_n\}$,
$$\lim_{n \to \infty} \sup_{s,t \in [0,1]}
\| g_n(e_n,t)g_n(e_n,s)\psi_V(g(s)) - \psi_V(g(s))g_n(e_n,t)g_n(e_n,s)\| = 0;$$
and $g$ is uniformly continuous, whence
$$\lim_{n \to \infty} \sup\{\|g(s)-g(t)\| : s,t \in [0,1], |s-t| \le
2/n\} = 0.$$
It follows that
\begin{eqnarray*}
\widetilde{W}_n(f) \widetilde{W}_n(g)
&=& \int_0^1\int_{(t-2/n) \vee 0}^{(t+2/n)\wedge 1}
g_n(e_n,t)\psi_V(f(t)g(t))g_n(e_n,t)g_n(e_n,s)^2 \, ds
\, dt + \mathrm{o}(1/n)\\
&=& \int_0^1\int_{0}^{1}
g_n(e_n,t)\psi_V(f(t)g(t))g_n(e_n,t)g_n(e_n,s)^2 \, ds
\, dt + \mathrm{o}(1/n) \\
& = &
\widetilde{W}_n(fg) + \mathrm{o}(1/n),
\end{eqnarray*}
as desired, where $\{\mathrm{o}(1/n)\}$ is a sequence of elements
whose norm tend to zero as $n$ tends to infinity.

(v). Because $e_me_n = e_m$ we have $e_mg_n(e_n,t) = e_mg_n(1,t)$, and
hence,
\begin{eqnarray*}
e_m W_n(f)e_m  - e_m \psi_V(f(1)) e_m  & = & e_m \big(\int_{1/n}^1
g_n(1,t)\psi_V\big(f(t) - f(1)\big)g_n(1,t) dt
\big)e_m  \\  & = & e_m \big(\int_{1-1/n}^1
g_n(1,t)\psi_V\big(f(t) - f(1)\big)g_n(1,t) dt
\big)e_m
\end{eqnarray*}
when $n \ge 2$. The norm of the last expression is less
than $\sup_{t \in [1-1/n,1]} \|f(t)-f(1)\|$.
\end{proof}

\begin{lemma} \label{lm:iota}
Suppose $A$ has property \emph{($\dagger$)}. Then
there is a \sh{} $\iota \colon C_1 \otimes \cK \to
A_\omega$ such that
$$I_{A_\omega}\big(\iota(J \otimes \cK)\big) = I_{A_\omega}(J)
$$
for all closed two-sided ideals $J$ in $C_1 \subseteq A_\omega$.
\end{lemma}

\noindent It follows in particular that $I_{A_\omega}(\iota(x \otimes
e)) = I_{A_\omega}(x)$ for all $x \in C_1$ and for all non-zero $e \in
\cK$.

\begin{proof} The multiplier algebra $\cM(A)$ is properly infinite by
  the assumption that $A$ is stable. Thus $\cM(A)_\omega$
  is properly
  infinite because it
  contains $\cM(A)_\omega$, and so we can find
  a sequence $s_1,s_2, \dots$ of isometries in $\cM (A)_\omega
  \subseteq \cM(A_\omega)$ with
  orthogonal ranges. Define $\iota \colon C_1 \otimes \cK \to
  A_\omega$ by $\iota(c \otimes e_{kl}) = s_kcs_l^*$, for $c \in C_1$
  and  where
  $\{e_{kl}\}_{k,l=1}^\infty$ is a system of matrix units for $\cK$.

Let $J$ be a closed two-sided ideal in $C_1$. Take $x$ in $J$. Then
$$x = s_1^*(s_1xs_1^*)s_1 = s_1^*\iota(x \otimes e_{11}) s_1 \in
I_{A_\omega}\big(\iota(J \otimes \cK)\big).$$
This proves $I_{A_\omega}(J) \subseteq I_{A_\omega}\big(\iota(J
\otimes \cK)\big)$. To show the reverse inclusion, we use that $J
\otimes \cK$ is the closure of the linear span of elements of the form
$x \otimes e_{kl}$, where $x \in J$ and $k,l \in \N$. Now, $\iota(x
\otimes e_{kl}) = s_kxs_l^* \in I_{A_\omega}(J)$, when $x$ belongs to
  $J$, and this proves $I_{A_\omega}\big(\iota(J
\otimes \cK)\big) \subseteq I_{A_\omega}(J)$.
\end{proof}

\noindent Put $S_n = \iota \circ W_n \circ \rho \colon A \to
A_\omega$. Each $S_n$ is a \cpc{}, and $\{S_n\}_{n=1}^\infty$ is
asymptotically multiplicative, ie., $\lim_{n \to
  \infty} \|S_n(ab) - S_n(a)S_n(b)\| = 0$ for all $a,b \in A$ by
Lemma~\ref{lm:W_n}~(iv).

\begin{lemma} \label{lm:s_n}
Suppose that $A$ has property \emph{($\dagger$)} and that $A$ is
nuclear. Then there are isometries $s_n$
in $\cM (A)_\omega$ such that $S_n(a) = s_n^*as_n$ for all $a \in A$.
\end{lemma}

\begin{proof} Doing the proof backwards, use
  Lemma~\ref{lm:1-step-inner} to reduce the problem to showing that
  $S_n$ is approximately 1-step inner. As $A$ is strongly purely
  infinite by assumption, so is $A_\omega$, cf.\
  \cite[Proposition~5.12]{KirRor:pi2}. By
  \cite[Theorem~7.21]{KirRor:pi2}, every
approximately inner \cpc{} $A \to A_\omega$ is actually approximately
1-step inner, so it suffices to show that $S_n$ is approximately
inner.

To see that $S_n$ is approximately inner we use
Corollary~\ref{cor:appr_inner} whereby it suffices to check that
  $S_n(a)$ belongs to $I_{A_\omega}(a)$ for every $a \in A$. Fix $a
  \in A$, and put $J =
  I_A(a)$, $J_1 = I_{C_1}(J \cap C)$, and $J_2 = I_{A_\omega}(a)$, and
  note that $J_1 \subseteq J_2$. Now, $\rho_t(a) \in J$
  for all $t \in [0,1]$, whence $\psi_V(\rho_t(a)) \in \cM(C_1 \otimes
  \cK, J_1 \otimes \cK)$ by \eqref{eq:psi_V}. It therefore follows
  from Lemma~\ref{lm:W_n}~(iii) that $W_n(\rho(a))$ belongs to $J_1
  \otimes \cK$, and so $S_n(a) = \iota\big(W_n(\rho(a))\big)$ belongs
  to $\iota(J_1 \otimes \cK)$. Now, $I_{A_\omega}(\iota(J_1 \otimes
  \cK)) = I_{A_\omega}(J_1) \subseteq J_2$ by
  Lemma~\ref{lm:iota}.
\end{proof}

\begin{lemma} \label{lm:T_n}
Suppose that $A$ has property \emph{($\dagger$)} and that $A$ is
nuclear.
There is a sequence of 1-step inner \cpc s $T_n \colon A_\omega \to
A_\omega$ such that $(T_n \circ S_n)(a) \to a$ for all $a \in A$.
\end{lemma}

\begin{proof} To ease the notation we let $B$ denote $A_\omega$ (so
  that $S_n \colon A \to B$).
It suffices to show that for each finite subset $F$ of
  contractions in $A$ and for each $\ep >0$ there is a natural number
  $n_0$ and a
  sequence $\{g_n\}_{n = n_0}^\infty$ of contractions in $B$
  such that $\|g_n^*S_n(a)g_n-a\| < \ep$ for all $a \in F$ and all $n
  \ge n_0$. Upon
  making a small perturbation of the elements in $F$ we may
  assume that $A$ contains a positive contraction $e$ such that
  $ea=ae=a$ for all $a \in F$. Put $F_1 = F \cup \{e\}$.

Recall that $B_\infty$ denotes the quotient algebra
$\ell_\infty(B)/c_0(B)$, that $\pi_\infty \colon \ell_\infty(B) \to
B_\infty$ is the quotient mapping, and  that $A \subseteq A_\omega = B
\subseteq B_\infty$ (where the canonical inclusion mapping are
suppressed).

Let $S \colon A \to B_\infty$ be given by
$$S(a) = \pi_{\infty}(S_1(a),S_2(a), S_3(a), \dots), \qquad a
\in A.$$
Then $S$ is a \sh{}, cf.\ Lemma~\ref{lm:W_n}~(iv). We show first that
$S$ is ideal preserving in the sense that
\begin{eqnarray} \label{eq:S}
a \in I_{B_\infty}(S(a)), \qquad a \in A.
\end{eqnarray}
It suffices to verify \eqref{eq:S} for elements $c$ in
$C$ (because $C$ separates ideals in $A$), so take such an element.
Retain the notation $\{e_m\}$ for the
approximate unit for $C_1 \otimes \cK$ that is used to define
$W_m$. Then
\begin{equation} \label{eq:r1}
\|e_{n-1}W_n(\rho(c))e_{n-1} - e_{n-1}\psi_V(c)e_{n-1}\| \le \sup \{
\|\rho_t(c) - \rho_1(c)\| : t \in [1 - 1/n,1] \}
\end{equation}
for all $n \ge 2$ by Lemma~\ref{lm:W_n}~(v). By \eqref{eq:Kasparov}, and since
$e_{n-1}$ commutes with $1 \otimes e_{11}$, we get
\begin{equation*} 
(1 \otimes e_{11})e_{n-1}\psi_V(c)e_{n-1}(1 \otimes e_{11}) =
e_{n-1}(V(c) \otimes e_{11})e_{n-1} = e_{n-1}(c \otimes e_{11})e_{n-1}
\underset{n \to \infty}{\to} c \otimes e_{11}.
\end{equation*}
Multiply \eqref{eq:r1} by $1 \otimes e_{11}$ from the left and the
right:
\begin{equation} \label{eq:r3}
\lim_{n \to \infty} \|(1 \otimes e_{11})e_{n-1} W_n(\rho(c))e_{n-1}(1
\otimes e_{11}) - c \otimes e_{11}\| = 0.
\end{equation}
Apply the map
$\iota \colon C_1 \otimes \cK \to A_\omega$ defined in
Lemma~\ref{lm:iota} to \eqref{eq:r3}, and recall from the proof of
that lemma that $\iota(c \otimes e_{11}) =
s_1cs_1^*$ for some isometry $s_1$ in $\cM(B)$,
\begin{eqnarray*}
\lim_{n \to \infty} \| x_n^*S_n(c)x_n - s_1cs_1^*\| = 0, \quad
\text{when} \; x_n = \iota(e_{n-1}(1 \otimes e_{11})) \in B.
\end{eqnarray*}
Hence $s_1^*x^*S(c)xs_1 = c$ when $x = \pi_\infty(x_1,x_2, \dots)$.
This proves \eqref{eq:S}.

Let $\sK$ be the set of \cp{} maps $A \to B_\infty$
of the form $a \mapsto d^*S(a)d$ for some $d$ in
$B_\infty$. We show next that $\sK$ is operator convex, cf.\
Definition~\ref{def:convex}. It is trivial that
(ii) in Definition~\ref{def:convex} holds and that $\alpha W$ belongs
to $\sK$ whenever $\alpha \in \R^+$ and $W \in \sK$.

By Lemma~\ref{lm:O_2-in-commutant} there is a unital embedding of
$\cO_\infty$ into $\cM(B_\infty) \cap
S(A)'$, and so there are isometries $t_1,t_2$ in  $\cM(B_\infty) \cap
S(A)'$ with $t_i^*t_j = 0$ when $i \ne j$. Take $W_1,
W_2 \in \sK$. Then $W_j(a) = d_j^*S(a)d_j$ for some $d_1,d_2$ in
$B_\infty$. Put $d = t_1d_1+t_2d_2 \in B_\infty$. Now, $W_1+W_2$
belongs to $\sK$ (so (i) in Definition~\ref{def:convex} holds) because
$$d^*S(a)d = \sum_{i,j=1}^2 d_i^*t_i^*S(a)t_jd_j =
\sum_{i,j=1}^2 d_i^*t_i^* t_j S(a)d_j = W_1(a) + W_2(a).$$

To verify Definition~\ref{def:convex}~(iii), let $W \in \sK$, $c_1,
\dots, c_n \in A$, and $b_1, \dots, b_n \in
B_\infty$ be given. Take $d$ in
$B_\infty$ such that $W(a) = d^*S(a)d$, and put $d_1
= \sum_{j=1}^n S(c_j)db_j$. Then
$$d_1^*S(a)d_1 =\sum_{i,j=1}^n b_i^*d^*S(c_i^*)S(a)S(c_j)db_j =
\sum_{i,j=1}^n b_i^*d^*S(c_i^*ac_j)db_j =  \sum_{i,j=1}^n
b_i^*W(c_i^*ac_j)b_j,$$
whence the mapping $a \mapsto 
\sum_{i,j=1}^n b_i^*W(c_i^*ac_j)b_j$ belongs to $\sK$.

It now follows from  Proposition~\ref{prop:HB} that the inclusion
mapping $j \colon A \to B_\infty$ belongs to the point-norm closure of
$\sK$, because
$$j(a) = a \in I_{B_\infty}(S(a)) = I_{B_\infty}\big(\{W(a) \mid W \in \sK
\}\big), \qquad a \in A,$$
by \eqref{eq:S}. Hence there exists an element $d$ in
$B_\infty$ such that $\|d^*S(a)d -a\| < \ep/3$ for all $a$ in the given
finite subset $F_1$ of $A$.

We show next that there is a \emph{contraction} $g$ in
$B_\infty$
with $\|g^*S(a)g -a\| < \ep $ for all $a$ in the finite subset $F$ of
$A$. Note first that $e$ belongs to $F_1$ and that
$\|d^*S(e)d-e\|
< \ep/3$. Put $f = S(e)^{1/2}d$ and put $g = \|f\|^{-1}f$. Then $f^*S(a)f =
d^*S(a)d$ for all $a \in F$ (because $ae=ea=a$ for all $a \in F$); and
$$\|f-g\| \; \le \; |\|f\|-1| \; \le \; |\|f\|^2-1| \; = \;
|\|d^*S(e)d\|-1| \; \le \; \|d^*S(e)d - e\| \; < \; \ep/3.$$
It follows that
\begin{eqnarray*}
\|g^*S(a)g -a\| & \le & \|g^*S(a)g-f^*S(a)f\|+\|f^*S(a)f-a\| \\
& \le & 2\|f-g\| + \|f^*S(a)f-a\| \; < \; \ep
\end{eqnarray*}
for all $a \in F$.

Write $g = \pi_\infty(g_1,g_2, \dots)$, where $g_1,g_2, \dots$
are contractions in $B$. Then
$$\limsup_{n\to \infty}\|g_n^*S_n(a)g_n - a\| = \|g^*S(a)g-a\| < \ep$$
for all $a \in F$. This completes the proof, cf.\ the remarks in the
first paragraph of the proof. \end{proof}

\begin{lemma} \label{lm:t_n}
Suppose that $A$ has property \emph{($\dagger$)} and that $A$ is
nuclear. Then
for each $n$ there is an isometry $t_n$ in $\cM(A)_\omega$
such that
$t_n^*S_n(a)t_n = (T_n \circ S_n)(a)$ for all $a$ in $A$.
\end{lemma}

\begin{proof} Fix $n$ and let $B$ be the separable sub-\Cs{} of
  $A_\omega$ generated by $S_n(A)$. Let $T'_n \colon B \to A_\omega$
  be the restriction of the 1-step inner \cpc{} $T_n \colon A_\omega
  \to A_\omega$ from Lemma~\ref{lm:T_n}. Then $T'_n$ is still a 1-step
  inner \cpc, so by Lemma~\ref{lm:1-step-inner} there is an isometry
  $t_n$ in $\cM(A)_\omega$
  such that $T'_n(b) = t_n^*bt_n$ for all $b \in B$.
\end{proof}

\begin{lemma} \label{lm:u_n}
Suppose that $A$ has property \emph{($\dagger$)} and that $A$ is nuclear.
Then there is a sequence $\{u_n\}_{n=1}^\infty$ of unitary elements in
$\cM(A)_\omega$
that satisfies $\|u_n^*S_n(a)u_n-a\| \to 0$ for all $a \in A$.
\end{lemma}

\begin{proof}
The conditions on $A$ imply that $A$ is $\cO_2$-absorbing as explained
below Lemma~\ref{lm:psi_V}.

Combining Lemmas~\ref{lm:s_n} and \ref{lm:t_n} we obtain sequences
$\{s_n\}_{n=1}^\infty$ and $\{t_n\}_{n=1}^\infty$ of isometries in
$\cM(A)_\omega$
such that $s_n^*as_n = S_n(a)$ and
$\|t_n^*s_n^*as_nt_n-a\| \to 0$ for all $a \in A$. The existence of
the unitaries $u_n \in \cM(A)_\omega$ follows from
Remark~\ref{rem:cantor-bernstein} and the fact established in
Lemma~\ref{lm:W_n}~(iv) that  $\{S_n\}$ is asymptotically
multiplicative.
\end{proof}

\noindent Consult \cite[Definition~5.1]{KirRor:pi2} for a formal
definition being ``strongly purely infinite'' (or see also the
comments above Proposition~\ref{prop:extension}).

\begin{theorem} \label{thm:pi=ah_0}
Let $A$ be a separable, nuclear, strongly purely infinite
\Cs{}. Suppose that $A$ is zero-homotopic in an ideal-system preserving way
(cf.\ Definition~\ref{def:homotopy-trivial}). Then $A$ is stable and
$\cO_2$-absorbing, i.e, $A \cong A \otimes
\cO_2 \otimes \cK$, and $A$ is the inductive limit of a sequence
\begin{equation} \label{eq:ah0}
\xymatrix{C_0((\Gamma_1,v_1),M_{k_1}) \ar[r]^{\varphi_1} &
  C_0((\Gamma_2,v_2),M_{k_2}) \ar[r]^{\varphi_2}  & C_0((\Gamma_3,v_3),M_{k_3})
      \ar[r]^-{\varphi_3} &  \cdots \ar[r] & A,}
\end{equation}
where each $(\Gamma_j,v_j)$ is a pointed graph, each $k_j$ is a
natural number, and each $\varphi_j$ is an injective \sh.

In particular, $A$ is an AH$_0$-algebra.
\end{theorem}

\begin{proof} It is shown below Lemma~\ref{lm:psi_V} that the given
  conditions on $A$ imply that $A$ is stable and $\cO_2$-absorbing.

To prove that $A$ is an inductive limit as stipulated, we use
Theorem~\ref{thm:ah_0} whereby it suffices to show that condition (iv)
of that theorem holds. Let $F = \{a_1,a_2, \dots, a_m\}$ be a finite
subset of $A$, let $\ep
>0$, and let $\omega$ be a free filter on $\N$. Let $S_n \colon A \to
A_\omega$ and $u_n \in \cM(A)_\omega$
be as constructed earlier in this section. Choose
$n$ large enough so that $\|u_n^*S_n(a)u_n-a\| < \ep/2$ for all $a \in
F$. Realize $C_1 \otimes \cK$ as the closure of $\bigcup_{k=1}^\infty C_1
\otimes M_k$, and set $B_{n,k} = u_n^*(\iota(C_1 \otimes M_k))u_n$
(cf.\ \eqref{eq:diagram}). Then $B_{n,k}$ is a sub-\Cs{} of $A_\omega$
isomorphic to $C_1 \otimes M_k$ and hence to $M_k(C_0(X))$, where $X$
is the spectrum of the Abelian \Cs{} $C_1$. As the closure of
$\bigcup_{k=1}^\infty B_{n,k}$  equals the image of $\Ad {u_n} \circ
\iota$ and hence contains $u_n^*S_n(a)u_n$ for all $a \in A$, we can
find a natural number $k$ and elements $b_1,b_2, \dots, b_m$ in
$B_{n,k}$ such that $\|u_n^*S_n(a_j)u_n-b_j\| < \ep/2$ for all
$j$. We now clearly have $\|a_j-b_j\| < \ep$, so property (iv) of
Theorem~\ref{thm:ah_0} holds.
\end{proof}

\section{Applications} \label{sec:applications}

\noindent We remind the reader of the construction of a \Cs{} $\A$
considered in \cite{Ror:piah}. Take a dense
sequence $\{t_n\}_{n=1}^\infty$ in $[0,1)$, and let $\A$ be the
inductive limit of the sequence
$$
\xymatrix{C_0([0,1), M_2) \ar[r]^-{\varphi_1} &  C_0([0,1), M_4)
  \ar[r]^-{\varphi_2} &  C_0([0,1), M_8) \ar[r]^-{\varphi_3} &  \cdots
  \ar[r] & \A,}
$$
where $\varphi_n(f)(t) = \diag(f(t), f(t \vee t_n))$. Put $A_n =
C_0([0,1), M_{2^n})$, and let $\varphi_{\infty,n} \colon A_n \to \A$
denote the inductive limit map.

It was shown in \cite{Ror:piah} that $\A$ is stable and that $\A
\cong \A \otimes \cO_\infty$. It is clear that $\A$ is separable and
nuclear. Note also that $\A$ by its definition is an inductive limit as
in the \eqref{eq:ah0} in Theorem~\ref{thm:pi=ah_0}.

For each $n \in \N$ and for each $t
\in [0,1]$ put $I^{(n)}_t = C_0([0,t),M_{2^n}) \subseteq A_n$, i.e.,
a function $f \in A_n$ belongs to $I^{(n)}_t$ if and only if
$f|_{[t,1]} \equiv 0$. As
$\varphi_n^{-1}(I^{(n+1)}_t) = I^{(n)}_t$ for all $n$ and $t$, there
is for each $t \in [0,1]$ a (unique) closed two-sided
ideal $I_t$ in $\A$ with
$\varphi_{\infty,n}^{-1}(I_t) =I^{(n)}_t$ for all $n$. It is shown in
\cite{Ror:piah} that any closed two-sided ideal in $\A$ is equal to
$I_t$ for some $t \in [0,1]$.

\begin{proposition} \label{prop:zero-homotopic}
The \Cs{} $\A$ is zero-homotopic in an ideal-system preserving way. In
particular, $\A \otimes \cO_2 \cong \A$.
\end{proposition}

\begin{proof} For each $s \in [0,1]$ let $\psi_n^{(s)} \colon A_n \to
  A_n$ be the
\sh{} given by $\psi_n^{(s)}(f)(t) = f(t \vee s)$. Note that
$\psi_n^{(s)}$ is a \sh{} for each $s$, that $s \mapsto
\psi_n^{(s)}(f)$ is continuous for all $f \in A_n$, that
$\psi_n^{(s)}(I_t \cap
A_n) \subseteq I_t \cap A_n$, and that $\psi_n^{(0)} = \Id$ and
$\psi_n^{(1)} = 0$. Moreover, as $t \vee t_n \vee s = t \vee s
\vee t_n$ for each $n$, we get a commutative diagram
$$
\xymatrix@C+1pc@R+.3pc{A_1 \ar[r]^-{\varphi_1} \ar[d]_-{\psi_1^{(s)}} & A_2
  \ar[r]^-{\varphi_2} \ar[d]_-{\psi_2^{(s)}} & A_3 \ar[r]^-{\varphi_3}
  \ar[d]_-{\psi_3^{(s)}} &  \cdots \ar[r] & \A \ar@{-->}[d]^-{\psi^{(s)}} \\
A_1 \ar[r]_-{\varphi_1}  & A_2
  \ar[r]_-{\varphi_2}  & A_3 \ar[r]_-{\varphi_3}
   &  \cdots \ar[r] & \A }
$$
that induces a point-norm continuous path $\psi^{(s)}$, $s \in [0,1]$,
of $^*$-endomorphisms on $\A$ that satisfies $\psi^{(s)}(I_t)
\subseteq I_t$ for all $s$ and $t$, and $\psi^{(0)} = \Id$ and
$\psi^{(1)} = 0$. This proves that $\A$ is zero-homotopic in an
ideal-system preserving way. Now Theorem~\ref{thm:Kir} yields that $A
\cong A \otimes \cO_2$.
\end{proof}

\begin{proposition} \label{prop:A_I-B}
The tensor product $\A \otimes B$
is zero-homotopic in an ideal-system preserving way for any separable
nuclear \Cs{} $B$. Hence the conclusion of
Theorem~\ref{thm:pi=ah_0} holds for $\A \otimes B$,
and in particular, $\A \otimes B$ is an AH$_0$-algebra as in \eqref{eq:ah0}.
\end{proposition}

\begin{proof} Let $\rho_t \colon \A \to \A$ be a point-norm continuous
  path of ideal preserving \sh s such that $\rho_1 = \Id$ and $\rho_0
  = 0$. We assert that $\sigma_t = \rho_t \otimes \Id_B$ then realizes
  an ideal-system preserving zero-homotopy of $\A \otimes B$. Clearly,
  $\sigma_1 = \Id$ and $\sigma_0 =0$. We must
  show that $\sigma_t(I) \subseteq I$ for
  every closed two-sided ideal $I$ in $\A \otimes B$ and for every
  $t$. If $x \otimes y \in I$ (where $x \in \A$ and $y \in B$), then
  $\rho_t(x)$ 
  belongs to the ideal in $\A$ generated by $x$, so 
  $\sigma_t(x \otimes y) = \rho_t(x) \otimes y$ belongs to $I$. This
  proves the claim, because 
   $I$ is the closure of the linear span of
  elementary tensors by Blackadar's \cite[Theorem~3.3]{Bla:tensor}.
\end{proof}

\begin{proposition} \label{prop:auto}
For each homeomorphism $h$ on $[0,1]$ with $h(1)=1$ there is an
automorphism $\alpha_h$ on $\A$ with $\alpha_h(I_t) = I_{h(t)}$.
\end{proposition}

\begin{proof}
Retain the notation set up in the construction of the \Cs{} $\A$ at
the beginning of this section.
For each $n$ consider the \sh{} $\beta_{h,n} \colon A_n=
C_0([0,1),M_{2^n}) \to \A$ given by $\beta_{h,n}(f) =
\varphi_{\infty,n}(f \circ h^{-1})$. 
Fix $f \in A_n$ and set $t = \inf\{s \in [0,1] \mid f|_{[s,1]} \equiv
0\}$. Then $I_{\A}(\varphi_{\infty,n}(f)) = I_t$ and $I_{\A}(\beta_{h,n}(f)) =
I_{h(t)}$. It follows that the closed two-sided ideal in
$\A$ generated by $\beta_{h,n}(I^{(n)}_t)$ is $I_{h(t)}$.

The two \sh s $\beta_{h,n+1} \circ \varphi_n$
and $\beta_{h,n}$ are not equal, but they are approximately unitarily
equivalent, cf.\ Theorem~\ref{thm:O_2-uniqueness}, because
$\beta_{h,n}(f)$ and $(\beta_{h,n+1} \circ \varphi_n)(f)$ generate the
same closed two-sided ideal in $\A$ for all $f \in A_n$, namely
$I_{h(t)}$, where $t$ as above is $\inf\{s \in [0,1] \mid f|_{[s,1]}
\equiv 0\}$.

By a one-sided approximate intertwining (after Elliott), see for
example \cite[Theorem~1.10.14]{Lin:amenable}, there is an endomorphism
$\beta_h$ on $\A$ such that $\beta_h \circ \varphi_{\infty,n}$ is
approximately unitarily equivalent to $\beta_{h,n}$ for all $n \in
\N$. The three sets
$$\beta_h(I_t), \qquad (\beta_h \circ \varphi_{\infty,n})(I^{(n)}_t),
\qquad \beta_{h,n}(I^{(n)}_t)$$
generate the same closed two-sided ideal in $\A$, namely $I_{h(t)}$.

Consider now the two \sh s $\beta_h$ and $\beta_{h^{-1}}$ of $\A$. By
the argument above, the closed two-sided ideals in $\A$ generated by
$(\beta_h \circ \beta_{h^{-1}})(I_t)$ and
$(\beta_{h^{-1}} \circ \beta_h)(I_t)$, respectively,
are both equal to $I_t$. We can
therefore conclude from Theorem~\ref{thm:O_2-uniqueness} that the
three \sh s $\beta_h \circ \beta_{h^{-1}}, \beta_{h^{-1}} \circ
\beta_h, \Id_{\A}$ are approximately unitarily equivalent. By a new
approximate intertwining argument, cf.\
\cite[Corollary~2.3.4]{Ror:encyc}, we obtain an automorphism
$\alpha_h$ on $\A$ which is approximately unitarily equivalent to
$\beta_h$. As the closed two-sided ideals generated by $\beta_h(I_t)$
and $\alpha_h(I_t)$ are equal, we conclude that $\alpha_h(I_t) =
I_{h(t)}$ as desired.
\end{proof}

\begin{lemma} \label{lm:AxZ}
Let $h$ be a homeomorphism on $[0,1]$ with $h(1)=1$, and let
$\alpha_h$ be an automorphism on $\A$ that satisfies
$\alpha_h(I_t)=I_{h(t)}$ for all $t \in [0,1]$, cf.\
Proposition~\ref{prop:auto}.
Then the crossed product $\A \rtimes_{\alpha_h} \Z$ is simple if (and only
if) $h$ has no fixed points in the open interval $(0,1)$.
\end{lemma}

\begin{proof}
  The assumptions imply that $\A$ is
  $(\Z,\alpha_h)$-simple (i.e., that there are no non-trivial ideals
  invariant under all powers of $\alpha_h$). It therefore follows from
  \cite[Theorem~7.2]{OlePed:C*-dynamicIII}
  that $\A \rtimes_{\alpha_h} \Z$ is simple if
  $\alpha_h^n$ is properly outer for all $n \ge 1$. Note that
  $\alpha_h^n = \alpha_{h^n}$ and that $h^n = h \circ h \circ \cdots
  \circ h$ has no fixed points in $(0,1)$
  when $n \ne 0$. Indeed,
  $h$ must be a strictly increasing continuous function
  on $[0,1]$ fixing $0$ and $1$ and
  satisfying either $h(t)<t$
  for all $t\in (0,1)$ or $h(t)>t$ for all
  $t\in (0,1)$. Hence $\A$ is $(\Z,\alpha_h^n)$-simple when $n \ne
  0$. It is a consequence of
  \cite[Theorem~6.6]{OlePed:C*-dynamicIII} that an automorphism,
  with no non-trivial invariant ideals, is properly outer if it is not
  approximately inner. As an approximately inner automorphism
  leaves all ideals invariant, we can conclude that $\alpha_h^n$ is
  properly outer for all $n \ne 0$. This completes the proof.
\end{proof}

\begin{proposition} \label{lm:AxZII}
There is an automorphism $\alpha$ on $\A$ such that $\A \rtimes_\alpha
\Z$ is isomorphic to $\cO_2 \otimes \cK$.
\end{proposition}

\begin{proof} Put $B = \A \otimes \cO_2$, and recall from
  Proposition~\ref{prop:zero-homotopic} that $B \cong \A$. Set $h(t) = t^2$,
  and use Proposition~\ref{prop:auto} to find an automorphism
  $\alpha_h$ on $\A$ with
  $\alpha_h(I_t) = I_{h(t)}$. Consider the automorphism $\alpha =
  \alpha_h \otimes \Id_{\cO_2}$ on $B$. The crossed product $B
  \rtimes_\alpha \Z$ is
  nuclear and separable, and it absorbs $\cO_2$:
$$B \rtimes_\alpha \Z = (\A \otimes \cO_2) \rtimes_{\alpha_h \otimes
  \Id_{\cO_2}} \Z \cong (\A \rtimes_{\alpha_h} \Z) \otimes \cO_2.$$
Moreover, $B \rtimes_\alpha \Z$ is stable because $B$ is stable (see
for example \cite[Corollary~4.5]{HjeRor:stable}).

We now identify $\A$ and $B$. Then $\alpha(I_t)
  = I_{t^2}$ for all $t \in [0,1]$, and so $\A \rtimes_\alpha \Z$ is
  simple by Lemma~\ref{lm:AxZ}.

It finally follows from \cite[Theorem~3.8]{KirPhi:classI} that  $\A
\rtimes_\alpha \Z$ is isomorphic to $\cO_2 \otimes \cK$.
\end{proof}

\begin{theorem} \label{thm:B}
Let $B$ be a separable, nuclear, stable, $\cO_2$-absorbing \Cs. It
follows that there is an $\cO_2$-absorbing \Cs{} algebra $A$, which is
zero-homotopic in an ideal-system preserving way, and which is an
inductive limit
$$\xymatrix{C_0((\Gamma_1,v_1),M_{k_1}) \ar[r]^{\varphi_1} &
  C_0((\Gamma_2,v_2),M_{k_2}) \ar[r]^{\varphi_2}  &
  C_0((\Gamma_3,v_3),M_{k_3})
      \ar[r]^-{\varphi_3} &  \cdots \ar[r] & A,}$$
where every $(\Gamma_j,v_j)$ is a pointed graph, $k_j$ is a
natural number, and  the connecting map $\varphi_j$ is a
$^*$-monomorphism, and there is an automorphism $\beta$ on $A$ such
that $B \cong A \rtimes_\beta \Z$.
\end{theorem}

\begin{proof} Take an automorphism $\alpha$ on $\A$ as in
  Proposition~\ref{lm:AxZII}. Put $A = \A \otimes B$, and put $\beta =
  \alpha \otimes \Id_B$. Then $A$ is as claimed by
  Proposition~\ref{prop:A_I-B}, and
$$A \rtimes_\beta \Z \, = \, (\A \otimes B) \rtimes_{\alpha \otimes
  \Id_B} \Z \, \cong \, (\A \rtimes_\alpha \Z) \otimes B \, \cong \,
\cO_2 \otimes \cK \otimes B \, \cong \, B.$$
\end{proof}

\noindent
We shall now use the fact that $B\otimes \A$ is an AH$_0$-algebra
for every separable, nuclear, strongly purely infinite \Cs{} $B$
to prove a step in the direction of a topological characterization
of the primitive ideal space of separable nuclear \Cs s.
The following notion of a regular subalgebra is crucial
for this.

\begin{definition} \label{def:regular-subalgebra}
A \sCs{} $C$ of  a \Cs{} $B$ is called \emph{regular}
in $B$, if the following hold for all closed ideals $I,J$ of $B$:
\begin{enumerate}
\item $I\cap C=J\cap C$ implies $I=J$,
\item $(I+J)\cap C = (I\cap C)+(J\cap C)$.
\end{enumerate}
\end{definition}

\noindent If $C \subseteq B$ is regular, then the map from the ideal
lattice of $B$ into the ideal lattice of $C$ given by $I \mapsto I
\cap C$ is an injective lattice morphism (i.e., this map preserves the two
lattice operations $\vee$ and $\wedge$). 
Moreover, this map also
preserves least upper bounds and greatest lower bounds of
\emph{infinite} families. Indeed, if $\{I_\alpha\}_{\alpha \in
  {\mathbb{A}}}$ is a family of ideals in $B$, then its greatest lower
bound is $\bigcap_\alpha I_\alpha$, and $\big(\bigcap_\alpha
I_\alpha\big) \cap C = \bigcap_\alpha (I_\alpha \cap C)$. If  
$\{I_\alpha\}_{\alpha \in {\mathbb{A}}}$ is upwards directed, 
then its least upper bound is $\overline{\bigcup_\alpha I_\alpha}$,
and
$$\Big(\overline{\bigcup_\alpha I_\alpha}\Big) \cap C =
\overline{\bigcup_\alpha I_\alpha \cap C}.$$
To form the least upper bound of a general family  $\{I_\alpha\}_{\alpha \in
  {\mathbb{A}}}$ one proceeds as above with the upper directed family
consisting of all finite sums of the ideals from the family $\{I_\alpha\}_{\alpha \in
  {\mathbb{A}}}$; and this operation is also preserved by our lattice
morphism by the results above.

In the case where $C$ is an Abelian regular sub-\Cs{} of $B$, so that
$C$ is isomorphic to $C_0(X)$ for some locally compact Hausdorff space
$X$, then, upon identifying the ideal lattice of $C$ with the lattice,
$\cO(X)$, of open subsets of $X$, we get an injective
lattice morphism from the ideal lattice of $B$ into $\cO(X)$.
This morphism maps $B$ to $X$ and $0$ to $\emptyset$. 
This morphism
also preserves least upper bounds and greatest lower bounds of infinite
families (the greatest lower bound of a family $\{U_\alpha\}_{\alpha \in
  \mathbb{A}}$ of open sets is the interior of $\bigcap_\alpha U_\alpha$).

Note that $(I\cap C)+(J\cap C) \subseteq (I+J)\cap C$ holds for any
sub-\Cs{} $C \subseteq B$ and for all closed two-sided ideals $I,J$ in
$B$. This inclusion can be strict, even when (i) above is
satisfied. 

We need some elementary lemmas
to get our final result, Theorem~\ref{thm:C-exist-for-nuc}.

\begin{lemma} \label{lm:permanence-of-regular}
Suppose that $A$, $B$, $D$ are \Cs{}s, such that $A \subseteq D$ and $D$
is simple and nuclear.
\begin{enumerate}
\item If $C_1\subseteq C_2\subseteq B$ and $C_1$ is
regular in $C_2$, then  $C_2$ is regular in $B$
if and only if $C_1$ is regular in $B$.
In particular, if $C_1$ is a full hereditary \sCs{} of $C_2$,
then $C_1$ is regular in $B$ if and only if
$C_2$ is regular in $B$.
\item
If $C_k\subseteq B_k$ are regular in $B_k$ for
$k=1,2$, then $C_1\oplus C_2$ is regular
in $B_1\oplus B_2$.
\item
If $C\subseteq B$ is regular in $B$,
then $C\otimes D$ is regular in $B\otimes D$.
\item If $F$ is a \Cs{} with
$B\otimes A\subseteq F \subseteq B \otimes D$ and $B\otimes A$ is full
in $F$, then $F$ is regular in $B\otimes D$.
\item
If $C$ is a regular \sCs{} of $B\otimes A$, then
$C$ is regular in $B\otimes D$.
\item If $E$ is the closure of the union of a sequence of \Cs s
$E_1\subseteq E_2\subseteq \ldots$, and if $C$ is a sub-\Cs{} of $E$
for which $E_n\cap C$ is regular
in $E_n$ for every $n\in \N$, and $\bigcup_n (E_n \cap C)$ is dense in
$C$, then $C$ is regular in $E$.
\end{enumerate}
\end{lemma}

\begin{proof} 
(i). To see that $C_2$ is regular in $B$ if $C_1$ is regular in $B$, take
two ideals $I,J$ of $B$ and use regularity of $C_1 \subseteq C_2$ and
of $C_1 \subseteq B$ to see that 
$$C_1 \cap (C_2 \cap I + C_2 \cap J) = C_1 \cap I + C_1 \cap J = C_1
\cap (I+J) = C_1 \cap(C_2 \cap (I+J)),$$
whence $C_2 \cap I + C_2 \cap J = C_2 \cap (I+J)$ again by regularity
of $C_1 \subseteq C_2$. The remaining claims of (i) and (ii) are
straightforward.  

In (iii), (iv), and (v) we need the following observation.
Since $D$ is exact and simple, every closed two-sided ideal
of $B\otimes D$ is of the form $I\otimes D$ for some closed two-sided
ideal $I$ of $B$. This was first proved by Blackadar,
\cite[Theorem 3.3]{Bla:tensor}, in the
case where $D$ is nuclear (all we need here).

(iii). This follows easily from the characterization of ideals of $B
\otimes D$  derived above and from the
identities: $I \otimes D + J \otimes D = (I+J) \otimes D$ and $(I \otimes
D) \cap (C \otimes D) = (I \cap C) \otimes D$, that hold for all closed
two-sided ideals $I$ and $J$ in $B$.

(iv). Property (i) of Definition~\ref{def:regular-subalgebra}
holds for $B \otimes A \subseteq B \otimes D$,
and hence also when $B \otimes A$ is replaced by the larger algebra
$F$. 

Let $F_1$ be the (full) hereditary sub-\Cs{} of $F$ generated by $B \otimes
A$. We show below that $F_1$ is regular in $B \otimes D$, and by (i)
this will entail that $F$ is regular in $B \otimes D$. 

We note first that 
$$(I \otimes D) \cap F_1 = \overline{(I \otimes 1) F_1}$$
for every closed two-sided ideal $I$ of $B$ (with $I \otimes 1$ viewed
as a subalgebra of $B \otimes \cM(D)$). If $x$ belongs to $(I \otimes
D) \cap F_1$, then $x$ belongs to the closure of $(I \otimes 1)x$ which
is contained in $\overline{(I \otimes 1) F_1}$. In the other direction,
$(I \otimes 1) F_1 \subseteq I \otimes D$. By construction of $F_1$ we
have $(B \otimes 1)F_1 \subseteq F_1$, whence $(I \otimes 1)F_1
\subseteq F_1$. 

Let now $I, J$ be two closed two-sided ideals of $B$. Then
\begin{eqnarray*}
(I\otimes D + J \otimes D) \cap F_1 & = & \big((I+J) \otimes D\big) \cap F_1
\; = \; \overline{\big((I+J) \otimes 1\big)F_1} \\ 
&= & \overline{(I \otimes 1)F_1 + (J \otimes 1) F_1} 
\;  = \; \overline{(I \otimes 1)F_1} +
\overline{(J \otimes 1)F_1} \\
&= & (I \otimes D) \cap F_1 + (J \otimes D) \cap F_1,
\end{eqnarray*}
as desired.

(v). Combine (i) and (iv) (with $F=B\otimes A$).

(vi). By a standard property of inductive limits, $K$ is the
closure of $\bigcup_{n=1}^\infty K \cap E_n$ whenever $K$ is a closed
two-sided ideal of $E$ or of $C$. Let now $I, J$ be two closed two-sided
ideals in $E$. If $I \cap C = J \cap C$, then $I \cap C \cap E_n = J
\cap C \cap E_n$ and hence $I \cap E_n = J \cap E_n$ for all $n$
(because $C \cap E_n \subseteq E_n$ is regular), and this implies
$I = J$. Next, using again that $C \cap E_n \subseteq E_n$ is regular,
we get
$$\big(I \cap E_n +J \cap E_n\big) \cap C  \; = \; I \cap E_n \cap C +
J \cap E_n \cap C \; \subseteq \; I \cap C + J \cap C,$$
for all $n$, whence
$$(I+J) \cap C \; = \; \Big(\overline{\bigcup_{n=1}^\infty I \cap E_n +J \cap
  E_n}\Big) \cap C \; \subseteq \;  I \cap C + J \cap C.$$
\end{proof}

\begin{lemma}\label{lm:MASA-of-building-blocks}
Suppose that $A$ is a building block
for an AH$_0$-algebra, i.e., $A$ is isomorphic to $\bigoplus_{k=1}^r
C_0(X_k, M_{n_k})$
for some locally compact spaces
$X_k$ and natural numbers $n_k$, $k=1,\ldots, r$. Then every maximal
Abelian sub-\Cs{} (masa) of $A$ is regular in $A$.
\end{lemma}

\begin{proof} By Lemma~\ref{lm:permanence-of-regular}~(ii)
it suffices to consider the case
  where $A$ has a single summand, i.e., $A = C_0(X,M_n)$. If $C$ is a
  masa in $A$, then $C$ necessarily
contains the center $C_0(X,\C 1_n)$
of $C_0(X, M_{n})$. Thus $C_k$ is regular in $C_0(X_k, M_{n_k})$
by Lemma \ref{lm:permanence-of-regular}~(iv),
with $D=M_{n_k}$ and $A=\C$.
\end{proof}

\noindent In general masa's need not be regular, there are even masa's
in sub-homogeneous \Cs{}s that are not regular.

For the proof of the next lemma note that any quotient
of a standard building block $\bigoplus_{k=1}^r C_0(X_k, M_{n_k})$ is
isomorphic to $\bigoplus_{k=1}^r C_0(Y_k, M_{n_k})$ for some closed
(possibly empty) subsets $Y_k$ of $X_k$ (with the convention
$C_0(\emptyset,M_n) = 0$).

\begin{lemma} \label{lm:C-exists-for-AH0}
Every (separable) AH$_0$-algebra $A$ contains an Abelian
\sCs{} $C\subseteq A$ which is regular in $A$.
\end{lemma}

\begin{proof} An AH$_0$-algebra is by definition an inductive limit of
  a sequence of algebras $B_1 \to B_2 \to B_3 \to \cdots$ where each
  $B_k$ is a building block (as described in
  Lemma~\ref{lm:MASA-of-building-blocks}). Let $\mu_{\infty,n} \colon
  B_n \to A$ be the inductive limit map and put $A_n =
  \mu_{\infty,n}(B_n)$. Then $A_1 \subseteq A_2 \subseteq A_3 \subseteq
  \cdots \subseteq A$, $\bigcup_{n=1}^\infty A_n$ is dense in $A$, and
  each $A_n$, being isomorphic to a quotient of $B_n$, is a building
  block as in Lemma~\ref{lm:MASA-of-building-blocks}.

Find inductively maximal Abelian \sCs{}s $C_n$
of $A_n$ such that $C_n\subseteq C_{n+1}$ for all $n$.
The $C_n$ is regular in $A_n$
by Lemma \ref{lm:MASA-of-building-blocks}~(i).
The closure $C$ of the union $\bigcup_{n=1}^\infty C_n$ is
an Abelian \sCs{} of $A$, and $C \cap A_n = C_n$ because $C_n$ is
maximal Abelian, and hence $C$ is regular in $A$ by
Lemma \ref{lm:MASA-of-building-blocks}~(vi).
\end{proof}

\begin{theorem} \label{thm:C-exist-for-nuc}
For every separable, nuclear \Cs{} $B$ there
exist an Abelian regular \sCs{} $C$ of $B\otimes \cO _2$ such that the
maximal ideal space $\widehat{C}=\mathrm{Prim}(C)$ has dimension at
most one.
\end{theorem}

\begin{proof} The \Cs{} $B  \otimes \A $ is an
AH$_0$-algebra by Proposition \ref{prop:A_I-B}.
Therefore, by Lemma \ref{lm:C-exists-for-AH0},
$B \otimes \A$ contains a regular  Abelian \sCs{} $C_1$.
Use any embedding of $\A$ into $\cO _2$ to obtain inclusions $C_1
\subseteq B  \otimes \A \subseteq B \otimes\cO_2$. 
By Lemma~\ref{lm:permanence-of-regular}~(v),
$C_1$ is regular in $B \otimes\cO _2$.
Thus $C_1 \otimes \cO_2 \subseteq B \otimes\cO_2\otimes \cO _2$
is regular by Lemma~\ref{lm:permanence-of-regular}~(iii).
There exists an Abelian \Cs{} $C$ such that $C_1\otimes 1 \subseteq C
\subseteq C_1\otimes \cO _2$ and such that $\dim(\widehat{C}) \le 1$ by
Proposition~\ref{prop:perturbation1}
(applied to $C_1 +\C 1 \cong C(X)$). 
Observe that $C_1 \otimes 1$ is full in $C$ because $(C_1 \otimes
1)C(C_1 \otimes 1)$ is dense in $C$.
It follows from
Lemma~\ref{lm:permanence-of-regular}~(iv) (with $A=\C 1_{\cO _2}$)
that $C$ is regular in $C_1\otimes \cO _2$.
Therefore, $C$ is also regular in $B\otimes  \cO_2\otimes
\cO_2 \cong B \otimes \cO_2$
by Lemma~\ref{lm:permanence-of-regular}~(i).
\end{proof}

\noindent It follows from Theorem~\ref{thm:C-exist-for-nuc} that if
$B$ is a separable nuclear \Cs{}, then there is a locally compact
second countable Hausdorff space $X$ with $\dim(X) \le 1$ such that
the ideal lattice of $B$ is order isomorphic to some sub-lattice
$\calL$ of the lattice of open subsets of $X$ with $X, \emptyset \in \calL$;
and this lattice morphism also preserves least upper bounds and
greatest lower bounds of infinite families. 
Indeed, let $C$ be
as in Theorem~\ref{thm:C-exist-for-nuc} and put $X =
\widehat{C}$. The map $I \mapsto (I \otimes \cO_2) \cap C$ is an
injective lattice morphism from the lattice of ideals of $B$ into the
lattice of ideals of $C$, cf.\ the remark below
Definition~\ref{def:regular-subalgebra}, and the latter is isomorphic
to the lattice of open subsets of $X$.

\begin{remark}\label{rem:on-HarKir:pi}
Theorem~\ref{thm:C-exist-for-nuc} will in \cite{HarKir:pi} be used by
the first-named author
and H.~Harnisch
to give several topological
characterizations of primitive ideal spaces of
separable nuclear \Cs{}s $B\,$, and to derive
an almost functorial construction
of nuclear \Cs{}s $B$ from given $T_0$-spaces $X$
in the described class.
More precisely, a $T_0$-space $X$ is the primitive
ideal space of a separable nuclear \Cs{} $B$
if and only if $X$ is a second countable
point-wise complete $T_0$-space and there
is a continuous, pseudo-open and pseudo-epic
map $\psi$ from a (one-dimensional)
locally compact Polish space $Y$
into $X$.

The map $\psi$ defines a certain partial order relation
$R_\psi$ on $Y$ (the pseudo-graph of $\psi$),
which in turn defines a Hilbert $C_0(X)$--bimodule
$H$ such that $X$ is naturally isomorphic
to the primitive ideal space of the Cuntz-Krieger-Pimsner
algebra $\cO _H$.
\end{remark}

\bibliographystyle{amsplain}
\providecommand{\bysame}{\leavevmode\hbox to3em{\hrulefill}\thinspace}
\providecommand{\MR}{\relax\ifhmode\unskip\space\fi MR }
\providecommand{\MRhref}[2]{%
  \href{http://www.ams.org/mathscinet-getitem?mr=#1}{#2}
}
\providecommand{\href}[2]{#2}

\newpage

\vspace{.3cm}
\noindent{\sc Humboldt-Universit{\"a}t zu Berlin, Institut f{\"u}r Mathematik,
Unter den Linden 6, D-10099 Berlin, Germany}

\vspace{.2cm}
\noindent{\sl E-mail address:} {\tt kirchbrg@mathematik.hu-berlin.de}\\

\vspace{.5cm}

\noindent{\sc Department of Mathematics, University of Southern
  Denmark, Odense,
  Campusvej~55, 5230 Odense M, Denmark}

\vspace{.3cm}

\noindent{\sl E-mail address:} {\tt mikael@imada.sdu.dk}\\
\noindent{\sl Internet home page:}
{\tt www.imada.sdu.dk/$\,\widetilde{\;}$mikael/welcome} \\

\end{document}